\documentclass[runningheads]{llncs}
\usepackage[dvipsnames]{xcolor}
\usepackage[sort&compress,numbers]{natbib}
\usepackage[colorlinks=true, allcolors=blue]{hyperref}
\usepackage[inline]{enumitem}
\usepackage{graphicx,array,pifont,soul,booktabs,multicol,setspace,euscript,subfigure,verbatim,xspace,standalone,pgfplots,bm,pgf,tikz,comment,url,float,longtable,amsfonts,amsmath,amssymb}
\usepackage[ruled,vlined,linesnumbered]{algorithm2e}
\usepackage[short,c2,nocomma]{optidef}
\usepackage[nameinlink]{cleveref}
\usepackage[super]{nth}
\setlist[enumerate]{label=(\roman*.)}
\crefname{line}{Line}{Lines}
\crefname{lemma}{Lemma}{Lemmata}
\crefname{theorem}{Theorem}{Theorems}
\crefname{proposition}{Proposition}{Propositions}
\crefname{algorithm}{Algorithm}{Algorithms}
\crefname{equation}{Equation}{Equations}
\crefname{definition}{Definition}{Definition}
\crefname{claim}{Claim}{Claim}
\crefname{corollary}{Corollary}{Corollaries}
\crefname{remark}{Remark}{Remarks}
\crefname{example}{Example}{Examples}
\crefname{figure}{Figure}{Figures}
\crefname{section}{Section}{Sections}
\crefname{table}{Table}{Tables}
\pgfplotsset{compat=1.16}

\definecolor{lime}{HTML}{A6CE39}
\DeclareRobustCommand{\orcidicon}{%
	\begin{tikzpicture}
	\draw[lime, fill=lime] (0,0) 
	circle [radius=0.16] 
	node[white] {{\fontfamily{qag}\selectfont \tiny ID}};
	\draw[white, fill=white] (-0.0625,0.095) 
	circle [radius=0.007];
	\end{tikzpicture}
	\hspace{-2mm}
}

\foreach \x in {A, ..., Z}{%
	\expandafter\xdef\csname orcid\x\endcsname{\noexpand\href{https://orcid.org/\csname orcidauthor\x\endcsname}{\noexpand\orcidicon}}
}

\newcommand{\KSet}{\mathcal{K}}
\newcommand{\KSetConv}{\conv(\mathcal{K})}
\newcommand{\NESet}{\mathcal{N}}
\newcommand{\EQFormulation}{\mathcal{E}}

\newcommand{\AGT}{\emph{AGT}\xspace}
\newcommand{\conv}{\operatorname{conv}}

\newcommand{\bd}{\operatorname{bd}}

\renewcommand{\int}{\operatorname{int}}
\newcommand{\ext}{\operatorname{ext}}
\newcommand{\proj}{\operatorname{proj}}

\newcommand{\POA}{\emph{PoA}\xspace}
\newcommand{\POS}{\emph{PoS}\xspace}

\newcommand{\KPG}{\emph{KPG}\xspace}
\newcommand{\NFG}{\emph{NFG}\xspace}
\newcommand{\CFLD}{\emph{CFLD}\xspace}
\newcommand{\PNE}{\emph{PNE}\xspace}
\newcommand{\IEDS}{\emph{IEDS}\xspace}
\newcommand{\PNEs}{\emph{PNE}s\xspace}

\newcommand{\MIP}{\emph{MIP}\xspace}

\newcommand{\QIPG}{\emph{qIPG}\xspace}
\newcommand{\QIPGs}{\emph{qIPG}s\xspace}
\newcommand{\IPG}{\emph{IPG}\xspace}
\newcommand{\IPGs}{\emph{IPG}s\xspace}

\newcommand{\NPH}{$\mathcal{NP}$-hard\xspace}
\newcommand{\SigmaTwoP}{$\Sigma^p_2$-hard\xspace}
\newcommand{\SigmaTwoPC}{$\Sigma^p_2$-complete\xspace}
 
\newtoggle{ARXIV}
\toggletrue{ARXIV}

\iftoggle{ARXIV}{
    \usepackage[margin=1.2in]{geometry}
    \usepackage[justification=centering]{caption}
}{
    \bibpunct[, ]{(}{)}{,}{a}{}{,}%
    \OneAndAHalfSpacedXII 
    \EquationsNumberedThrough
    \TheoremsNumberedThrough 
    \MANUSCRIPTNO{000000}
    
}
\begin{document}

\newcommand{\orcidauthorA}{0000-0002-5447-6245}
\newcommand{\orcidauthorB}{0000-0001-6612-1890}
\newcommand{\theTitle}{The ZERO Regrets Algorithm: Optimizing over Pure Nash Equilibria via Integer Programming}
\newcommand{\theTitleRunning}{Optimizing over Pure Nash Equilibria via Integer Programming}
\newcommand{\theAbstract}{Designing efficient algorithms to compute Nash equilibria poses considerable challenges in Algorithmic Game Theory and Optimization. In this work, we employ integer programming techniques to compute Nash equilibria in Integer Programming Games, a class of simultaneous and non-cooperative games where each player solves a parametrized integer program. We introduce \emph{ZERO Regrets}, a general and efficient cutting plane algorithm to compute, enumerate, and select Nash equilibria. Our framework leverages the concept of \emph{equilibrium inequality}, an inequality valid for any Nash equilibrium, and the associated \emph{equilibrium} separation oracle. We evaluate our algorithmic framework on a wide range of practical and methodological problems from the literature, providing a solid benchmark against the existing approaches. 
}

\iftoggle{ARXIV}{%
\title{\theTitle}

\author{
Gabriele Dragotto and 
Rosario Scatamacchia
}
\institute{
\orcidA{} \texttt{gdragotto@princeton.edu}\\
\orcidB{} \texttt{rosario.scatamacchia@polito.it}
}
\titlerunning{\theTitleRunning}

\maketitle     
\begin{abstract}
\theAbstract
\end{abstract}      \bibliographystyle{abbrvnat}
}{
\ARTICLEAUTHORS{%
    \AUTHOR{Gabriele Dragotto \orcidA{}}
    \AFF{ORFE, Princeton University and DS4DM, Polytechnique Montr\'eal \EMAIL{gdragotto@princeton.edu}
    }
    \AUTHOR{Rosario Scatamacchia \orcidB{}}
    \AFF{DIGEP, Politecnico di Torino \EMAIL{rosario.scatamacchia@polito.it}
    }
} 

\ABSTRACT{\theAbstract}
\TITLE{\theTitle}
\RUNAUTHOR{Dragotto and Scatamacchia}
\RUNTITLE{\theTitleRunning}
\bibliographystyle{informs2014}

\maketitle
 }

\section{Introduction}
Several real-world problems often involve a series of selfish agents optimizing their benefits while mutually affecting their decisions. The concept of Nash equilibrium \citep{nash_equilibrium_1950,nash_non-cooperative_1951} revolutionized the understanding of the agents' strategic behavior by proposing a flexible and interpretable solution, with consequences and applications in many different contexts. The Nash equilibrium constitutes a \emph{stable} solution, meaning that no single agent has an incentive to defect from it profitably. Nash equilibria, however, may intrinsically differ in their features, for instance, in terms of a given welfare function measuring the common good for the collectivity of the agents. Above all, the quality of equilibria often does not match the quality of the \emph{social optimum}, i.e., the best possible solution for the collectivity. In general, the social optimum is not a stable solution and, therefore, does not emerge naturally from the agents' interactions. Nevertheless, in numerous contexts, a central authority may suggest solutions to the agents, preferably ensuring that such solutions satisfy two foremost properties. First, the authority should ensure that little to no incentives exist for the agents to refuse the proposed solution. Second, the solution should be sufficiently close -- in terms of quality -- to the social optimum. The best trade-off between these two properties is the \emph{best} Nash equilibrium, i.e., a solution that optimizes a welfare function among the equilibria. Often, the main focus is on selecting a Pure Nash Equilibrium (\PNE), a stable solution where each agent selects one alternative with probability one (in contrast to a Mixed-Strategy equilibrium, where agents randomize over the set of their alternatives).
In this context, the Algorithmic Game Theory (\AGT) community pioneered the study of the interplay between Game Theory and algorithms with a focus on equilibria's efficiency \citep{nisan_algorithmic_2008}. The discipline attracted significant attention from the computer science and optimization communities, especially to study games where agents solve optimization problems (e.g., \citet{facchinei_finite-dimensional_2004}).
Several recent works \citep{hemmecke_nash-equilibria_2009,koppe_rational_2011,sagratella_computing_2016,carvalho_computing_2017,cronert_equilibrium_2020,schwarze_branch-and-prune_2022,harks_generalized_2022,dragotto_zero_2021} considered Integer Programming Games (\IPGs), namely games where the agents solve parametrized integer programs. In this work, we focus on a class of \emph{simultaneous} and \emph{non-cooperative} \IPGs among $n$ players (agents), as described in \cref{def:Game}, where each player controls $m$ integer variables. %
\begin{definition}[\IPG]
Each player $i=1,2,\dots,n$ solves \eqref{eq:def:game}, where $u^i(x^i, x^{-i})$ -- given $x^{-i}$ -- is a function in $x^i$ with integer coefficients, $A^i \in \mathbb{Z}^{r\times m}$, $b^i \in \mathbb{Z}^r$.
\begin{align}
    \max_{x^i} \{ u^i(x^i, x^{-i}): x^i \in \mathcal{X}^i \}, \; \mathcal{X}^i := \{ A^i x^i \le b^i, x^i \in \mathbb{Z}^{m} \}.
    \label{eq:def:game}
\end{align}
\label{def:Game}
\end{definition}
\noindent As standard game-theory notation, let $x^i$ denote the vector of variables of player $i$, and let the operator $\left( \cdot \right)^{-i}$ be $\left( \cdot \right)$ except $i$. The vector $x^{-i}=(x^1, \dots, x^{i-1},$ $x^{i+1}, \dots x^n)$ represents the variables of $i$'s \emph{opponents} (all players but $i$), and the set of linear constraints $A^i x^i \le b^i$ defines the feasible region $\mathcal{X}^i$ of player $i$. We assume all integer variables are lower and upper bounded, and thus that  $\mathcal{X}^i$ is finite. In \IPGs, the strategic interaction occurs in the players' objective functions, and not within their feasible regions.  
Specifically, players choose their strategy simultaneously, and each player $i$'s utility (or payoff) $u^i(x^i, x^{-i})$ is a function in $x^i$ parametrized in $i$'s opponents variables $x^{-i}$. Without loss of generality, we assume the entries of $A^i$ and $b^i$ and the coefficients of $u^i(x^i, x^{-i})$ are integers. Further, considering the space of all players' variables $(x^1,\dots,x^n)$, we assume one can always linearize the non-linear terms in each $u^i$ with a finite number of inequalities and auxiliary variables (e.g., \citet{sherali_reformulation-linearization_1999,vielma_mixed_2015}). We remark that this assumption is not restrictive; on the contrary, it enables us to tackle several games where the players' utilities are not linear (see \cref{sec:games}). 
Besides, we assume \begin{enumerate*} \item players have \emph{complete information} about the structure of the game, i.e., each player knows the other players' optimization problems via their feasible regions and objectives, \item each player is \emph{rational}, namely it always selects the best possible strategy given the information available on its opponents, and \item \emph{common knowledge of rationality}, namely each player knows its opponents are rational, and there is complete information \end{enumerate*}. 
\IPGs extend traditional resource-allocation tasks and combinatorial optimization problems to a multi-agent setting, and their modeling power lies precisely in the discrete variables and game dynamics they can model. Indeed, in several real-world applications, requirements such as indivisible quantities and fixed production costs often require the use of discrete variables (see, for instance, \citet{bikhchandani_competitive_1997}). Several recent works explored the application of \IPGs in various contexts. To name a few, \citet{gabriel_solving_2013} modeled energy production games, 
\citet{david_fuller_alternative_2017} proposed discrete unit commitment problem with fixed production costs,
\citet{anderson_supplier_2017} modeled a game where firms reserve discrete blocks of capacities from their suppliers, \citet{federgruen_multi-product_2015} proposed a price competition framework with $n$ competitors offering a discrete number of substitutable products, and \citet{carvalho_nash_2017} exploited \IPGs in the context of kidney exchange programs.
Despite the high potential impact of \IPGs in many domains, practitioners and researchers often make restrictive assumptions about the game's structure to guarantee that solutions are unique or computationally tractable. This is mainly due to the lack of a general, scalable and reliable methodology to select efficient solutions in \IPGs, which could potentially open new opportunities in terms of applications.
This lack is the core motivation behind our work: providing a general-purpose algorithmic framework to optimize over the solutions of \IPGs. Specifically, we focus on optimizing over the set of \PNEs for the \IPGs defined above and on characterizing the polyhedral structure of the set containing the \PNEs. The algorithmic framework possesses a solid theoretical foundation, and it integrates with the existing tools from the theory and practice of integer programming and combinatorial optimization. From a computational perspective, it is highly flexible, and it generally outperforms the algorithms available in the present literature.
Our framework is problem-agnostic and general, yet, it can be customized to address problem-specific needs.

\paragraph{Literature. } \citet{koppe_rational_2011} pioneered \IPGs by laying down their first formal definition. The authors also provided an algorithmic framework to enumerate \PNEs when the players' utilities are differences of piecewise linear convex payoff functions. Although their approach is theoretically well-grounded, there is no computational evidence of its effectiveness.
Indeed, even in some 2-player games (e.g., normal-form \citep{rosenberg_enumeration_2005,avis_enumeration_2010} and bimatrix \citep{audet_enumeration_2006} games) there are considerable computational challenges involved in the design of efficient algorithms for computing and selecting equilibria.
\citet{sagratella_computing_2016} proposed a branching method to enumerate the \PNEs in \IPGs where each $u^i(x^i,x^{-i})$ is convex in $x^i$. More recently, \citet{schwarze_branch-and-prune_2022} extended the work of \citet{sagratella_computing_2016} by proposing an improved branch-and-prune scheme that also drops the convexity assumption on the players' utilities.  \citet{del_pia_totally_2017} focused on totally-unimodular congestion games, namely \IPGs where players have totally-unimodular constraint sets $\mathcal{X}^i$. They propose a strongly-polynomial time algorithm to find a \PNE and derive some computational complexity results. Their results have been extended by \citet{kleer_computation_2021}. More recently, \citet{Dragotto_2021_CNP} proposed a general-purpose cutting plane algorithm to compute a \PNE in \IPGs where each player utility is linear in their variables and bilinear with respect to the other players' variables. However, their approach does not handle equilibria selection, and requires a specific structure on the players' objectives to derive the Karush–Kuhn–Tucker conditions associated with the linear relaxation of their optimization problems.
An important family of techniques for computing Mixed-Strategy equilibria is the one of support enumeration algorithms. The core idea is to determine if an equilibrium with a given support for each player -- e.g., a subset of its strategies -- exists in a normal-form game by solving a linear system of inequalities. \citet{porter_simple_2008} and \citet{sandholm_mixed-integer_2005} exploited this idea in the context of $n$-players normal-form games. Since equilibria in such games tend to have small supports, as proved theoretically by \citet{mclennan_expected_2005}, support enumeration algorithms tend to be practically efficient in normal-form games.
Inspired by the approach of \citet{porter_simple_2008}, \citet{carvalho_computing_2017} introduced the sample generation method ($SGM$) to compute an equilibrium in separable \IPGs (i.e., where each player's payoff takes the form of a sum-of-products) where players have bounded strategy sets. Their algorithm iteratively refines a sample of players' supports to compute an equilibrium or a correlated equilibrium (i.e., a generalization of the Nash equilibrium). However, the $SGM$ does not handle the enumeration or selection of equilibria, nor can it prove that no equilibrium exists. \citet{cronert_equilibrium_2020} modified the $SGM$ -- extending the work of \citet{carvalho_computing_2017} -- by proposing an enumerative algorithm to compute all the equilibria with the additional assumptions that all the players' variables are integer. They further complemented their approach with some considerations stemming from the theory of equilibria selection of \citet{harsanyi_new_1995}. 
Nevertheless,  identifying the correct samples leading to equilibria in \IPGs could be computationally cumbersome. While our approach shares a few elements with \citet{cronert_equilibrium_2020}, it does not require any sampling in order to compute and select equilibria. This fundamental aspect leads to significant differences in terms of practical effectiveness and performance of the algorithms (see \cref{sec:games}).

Although the previous methodological works provide an insightful perspective on the computability and the selection of equilibria in \IPGs, there are other significant intrinsic questions concerning the general nature of equilibria. 
Indeed, from the \AGT standpoint, not all equilibria are created equal. Three paradigmatic questions in \AGT and Game Theory are often: \begin{enumerate*} \item Does at least one \PNE exist? \item How good (or bad) is a \PNE compared to the social optimum? \item If more than one equilibrium exists, can one select the best \PNE according to a given measure of quality?  \end{enumerate*} 
Establishing that a \PNE does not exist may turn out to be a difficult task \citep{daskalakis_complexity_2009}. Nash proved that there always exists a Mixed-Strategy equilibrium in finite games, i.e., games with a finite number of strategies and players.
In \IPGs, where the set of players' strategies is large, deciding if a \PNE exists is generally a \SigmaTwoP decision problem in the polynomial hierarchy \citep{vaz_existence_2018}.
To measure the efficiency of equilibria, \citet{goos_worst-case_1999} introduced the concept of \emph{Price of Anarchy} (\POA), the ratio between the welfare value of the worst-possible equilibrium and the welfare value of a social optimum. Similarly, \citet{schulz_performance_2003,anshelevich_near-optimal_2003} introduced the \emph{Price of Stability} (\POS), the ratio between the welfare value of the best-possible equilibrium and a social optimum's one.  In the \AGT literature, many works focus on providing theoretical bounds for the \POS and the \POA, often by exploiting the game's structural properties \citep{anshelevich_near-optimal_2003,anshelevich_price_2008,chen_network_2006,nisan_algorithmic_2008,roughgarden_bounding_2004}. However, in practice, one may be interested in establishing the exact values of such prices in order to characterize the efficiency of equilibria in specific applications. This further highlights the need for general and effective algorithmic frameworks to select equilibria.

\paragraph{Contributions. }In this work, we shed new light on the intersection between \AGT and integer programming. We propose a new theoretical and algorithmic framework to efficiently and reliably compute, enumerate, and select \PNEs for the \IPGs in \cref{def:Game}. 
We summarize our contributions as follows:
\begin{enumerate}
    \item From a theoretical perspective, we provide a polyhedral characterization of the convex hull of the \PNEs. We adapt the concepts of valid inequality, closure, and separation oracle to the domain of Nash equilibria. Specifically, we introduce the concept of \emph{equilibrium inequality} to guide the exploration of the set of \PNEs. With this respect, we provide a general class of equilibrium inequalities and prove -- through the concept of \emph{equilibrium closure} -- they are sufficient to define the convex hull of the \PNEs.
    From a game-theory standpoint, we explore the interplay between the concept of rationality and cutting planes through the \emph{equilibrium inequalities}. Since in any game, a player $i$ may never play some of its strategies due to their induced payoffs, it is reasonable to think that player $i$ would only pick its strategies from a \emph{rational} subset of $\mathcal{X}^i$. In other words, we provide an interpretable criterion -- in the form of a cutting plane --  for a player to play or not some strategies. In this sense, what we propose constitutes an analytical and geometrical characterization of the sets of equilibria providing a novel perspective on equilibria selection.
    \item From a practical perspective, we design a cutting plane algorithm -- \emph{ZERO Regrets} -- that computes the most efficient \PNE for a given welfare function. This algorithm is flexible and scalable, it can potentially enumerate \emph{all} the \PNEs and compute approximate \PNEs. The algorithm exploits an \emph{equilibrium separation oracle},  a procedure separating non-equilibrium strategies from \PNEs through general and problem-specific equilibrium inequalities. Furthermore, our framework smoothly integrates with existing mathematical programming solvers, allowing practitioners to exploit the capabilities of the available optimization technologies.
    \item We evaluate our algorithmic framework on a range of applications and problems from the relevant works in the literature. We provide a solid benchmark against the existing approaches and show the flexibility and effectiveness of \emph{ZERO Regrets}. The classes of games we select derive from practical applications (e.g., competitive facility locations, network design) and methodological studies and the associated benchmark instances (e.g., games among quadratic programs). First, we consider the Knapsack Game, an \IPG where each player solves a binary knapsack problem. For this problem, we also provide theoretical results on the computational complexity of establishing the existence of \PNEs and two problem-specific equilibrium inequalities. Second, we focus on a Network Formation Game, a well-known and intensely investigated problem in \AGT, where players build a network over a graph via a cost-sharing mechanism. Third, we consider a Competitive Facility Location and Design game, where several sellers strategically decide the location and design of their facilities in order to maximize their revenues. Finally, we test our algorithm on a game where players solve integer problems with convex and non-convex quadratic objectives. \emph{ZERO Regrets} outperforms any baseline, proving to be highly efficient in both enumerating and selecting \PNEs.

\end{enumerate}
We remark that our framework can be extended to the non-linear case, i.e., when $u^i$ is non-linearizable. However, we focus on the linear case 
\begin{enumerate*}
\item to provide geometrical, polyhedral, and combinatorial insights on the structure of Nash equilibria in \IPGs, and
\item to foster the interaction with existing streams of research in Combinatorial Optimization. 
\end{enumerate*}

We structure the paper as follows. In \cref{sec:definitions}, we introduce the fundamental definitions and terminology. In \cref{sec:eqInequalities} we introduce the theoretical elements of our algorithmic framework. In \cref{sec:CompPNE}, we describe our cutting plane algorithm and its separation oracle and their extensions to compute approximate equilibria. In \cref{sec:games} we present an extensive computational campaign on the applications mentioned above, and, in \cref{sec:concluding}, we provide some concluding remarks.

\section{Definitions}
\label{sec:definitions}
We assume the reader is familiar with basic concepts of polyhedral theory and integer programming \citep{conforti_integer_2014}. We introduce the notation and definitions related to an \IPG instance $G$, where we omit explicit references to $G$ when unnecessary. Let $\mathcal{X}^i$ be the \emph{set of feasible strategies} (or the feasible set) of player $i$, and let any strategy $\bar{x}^i \in \mathcal{X}^i$ be a \emph{(pure) strategy} for $i$.
Any $\bar{x}=(\bar{x}^1,\dots,\bar{x}^n)$ -- with $\bar{x}^i \in \mathcal{X}^i$ for any $i$ -- is a \emph{strategy profile}. 
Let the vector $x^{-i}=(x^1,\dots, x^{i-1}, x^{i+1}, \dots x^n)$ denote the vector of \emph{the $i$'s opponents (pure) strategies}. 
The \emph{payoff} for $i$ under the profile $\bar{x}$ is   $u^i(\bar{x}^i,\bar{x}^{-i})$. We define $S(\bar{x})=\sum_{i=1}^n u^i(\bar {x}^i,\bar{x}^{-i})$ as the \emph{social welfare} corresponding to a given strategy profile $\bar{x}$.

\paragraph{Equilibria and Prices.}
A strategy $\bar{x}^i$ is a best-response strategy for player $i$ given its opponents' strategies $\bar{x}^{-i}$ if $u^i(\bar{x}^i,\bar{x}^{-i}) \ge u^i(\hat{x}^i,\bar{x}^{-i})$ for any  $\hat{x}^i \in \mathcal{X}^i$; equivalently, we say $i$ cannot profitably deviate to any $\hat{x}^i$ from $\bar{x}^i$. The difference  $u^i(\bar{x}^i,\bar{x}^{-i}) - u^i(\hat{x}^i,\bar{x}^{-i})$ is called the \emph{regret} of strategy $\hat{x}^i$ under $\bar{x}^{-i}$.
Let $\mathcal{BR}(i,\bar{x}^{-i})=\{ x^i \in \mathcal{X}^i: u^i(x^i,\bar{x}^{-i})\ge u^i(\hat{x}^i,\bar{x}^{-i}) \; \forall \; \hat{x}^i \in \mathcal{X}^i \}$ be the set of best-responses for $i$ under $\bar{x}^{-i}$. A strategy profile $\bar{x}$ is a \PNE if, for any player $i$ and any strategy $\hat{x}^i \in \mathcal{X}^i$, $u^i(\bar{x}^i,\bar{x}^{-i}) \ge u^i(\hat{x}^i,\bar{x}^{-i})$, i.e. any $\bar{x}^i$ is a best-response to $\bar{x}^{-i}$ (all regrets are $0$). Equivalently, in a \PNE, no player $i$ can unilaterally improve its payoff by deviating from its strategy $\bar{x}^i$. We define the \emph{optimal social welfare} as $OSW=\max_{x^1,\dots,x^n}\{ S(x) : x^i \in \mathcal{X}^i \; \forall i=1,2,\dots,n\}$. Given $G$, we denote as $\NESet=\{x=(x^1,\dots,x^n): x \text{ is a \PNE for }  G\}$ the set of its \PNEs. Also, let $\NESet^i:=\{ x^i : (x^i,x^{-i}) \in \NESet \}$, with $\NESet^i \subseteq \mathcal{X}^i$ be the set of \emph{equilibrium strategies} for $i$, namely the strategies of $i$ appearing in at least a \PNE. If $\NESet$ is not empty, let: \begin{enumerate*} \item $\dot{x} \in \NESet$ be so that $S(\dot{x})\le S(\bar{x})$ for any $\bar{x} \in \NESet$ (i.e., the \PNE with the \emph{worst} welfare), and \item  $\ddot{x} \in \NESet$ be so that $S(\ddot{x})\ge S(\bar{x})$ for any $\bar{x} \in \NESet$ (i.e., the \PNE with the \emph{best} welfare) \end{enumerate*}. Assuming w.l.o.g. $OSW > 0$ and $S(\ddot{x}) > 0$, the \POA of $G$ is $\frac{OSW}{S(\dot{x})}$, and the \POS is
$\frac{OSW}{S(\ddot{x})}$. The definitions of \POA and \POS hold when agents maximize a welfare function. Otherwise, when agents minimize their costs (e.g., the costs of routing packets in a network), we exchange numerator and denominator in both the \POA and the \POS.

\paragraph{Polyhedral Theory.} For a set $S$, let $\conv(S)$ be its convex hull. Let $P$ be a polyhedron: $\bd(P)$, $\ext(P)$,  $\int(P)$, are the boundary, the set of vertices (extreme points), and the interior of $P$, respectively. Let $P \subseteq \mathbb{R}^p$ and $\tilde{x} \notin P$ a point in $\mathbb{R}^p$. A \emph{cut} is a valid inequality $\pi^\top x \le \pi_0$ for $P$  violated by $\tilde{x}$, i.e., $\pi^\top \tilde{x} > \pi_0$ and $\pi^\top x \le \pi_0$ for any $x \in P$. 
Given a point $\hat{x} \in \mathbb{R}^p$ and $P$, we define the \emph{separation problem} as the task of determining that either \begin{enumerate*} \item $\hat{x} \in P$, or \item $\hat{x} \notin P$ and returning a cut $\pi^\top x \le \pi_0$ for $P$ and $\hat{x}$ \end{enumerate*}. For each player $i$, the set $\conv(\mathcal{X}^i)$ is the \emph{perfect formulation} of $\mathcal{X}^i$, namely an integral polyhedron whose vertices are in $\mathcal{X}^i$. 

\section{Lifted Space and Equilibrium Inequalities}
\label{sec:eqInequalities}
Cutting plane methods are attractive tools for integer programs, both from a theoretical and an applied perspective. The essential idea is to iteratively refine a relaxation of the original problem by cutting off fractional solutions via valid inequalities for the integer program's perfect formulation.
Nevertheless, in an \IPG where the solution paradigm is the Nash equilibrium, we argue there exist stronger families of cuts, yet, not necessarily \emph{valid} for each player's perfect formulation $\conv(\mathcal{X}^i)$. 
In fact, for any player $i$, some of its best-responses in $\bd(\conv(\mathcal{X}^i))$ may never appear in a \PNE, since no equilibrium strategies $\mathcal{N}^{-i}$ of $i$'s opponents induce $i$ to play such best-responses. In this work, we introduce a general class of inequalities to characterize the nature of $\conv(\mathcal{N})$. Such inequalities play a pivotal role in the cutting plane algorithm of Section \ref{sec:CompPNE}. 

\paragraph{Dominance and Rationality. }  We ground our reasoning in the concepts of \emph{rationality} and \emph{dominance} \citep{ pearce_rationalizable_1984,bernheim_rationalizable_1984}.
Given two strategies $\bar{x}^i \in \mathcal{X}^i$ and $\hat{x}^i \in \mathcal{X}^i$ for player $i$, $\bar{x}^i$ is strictly \emph{dominated} by $\hat{x}$ if, for any choice of opponents strategies $x^{-i}$, then $u^i(\hat{x},x^{-i}) > u^i(\bar{x},x^{-i})$. Then, a \emph{rational} player will never play dominated strategies. This also implies no player $i$ would play any strategy in $\int(\conv(\mathcal{X}^i))$. Since dominated strategies -- by definition  -- are never best-responses, they will never be part of any \PNE. In \cref{ex:KPGBoundDominance}, the set $\mathcal{X}^2$ is made of $3$ strategies $(x^2_1,x^2_2)=(0,0)$, $(1,0)$, $(0,1)$. 
Yet, $(x^2_1,x^2_2)=(0,0)$ is dominated by $(x^2_1,x^2_2)=(0,1)$, and the latter is dominated by $(x^2_1,x^2_2)=(1,0)$. However, when considering player 1, we need the assumption of \emph{common knowledge of rationality} to conclude which strategy the player will play. Player 1 needs to know that player 2 would never play $x^2_2 = 1$ to declare $(x^1_1,x^1_2)=(0,1)$ being dominated by $(x^1_1,x^1_2)=(1,0)$. When searching for a \PNE in this example, it follows that $\mathcal{N}^1=\{(x^1_1,x^1_2)=(1,0)\}$ and $\mathcal{N}^2=\{(x^2_1,x^2_2)=(1,0)\}$. This inductive (and iterative) process of removal of strictly dominated strategies is known as the \emph{iterated elimination of dominated strategies} (\IEDS). This process produces tighter sets of strategies and never excludes any \PNE from the game \citep[Ch.4]{tadelis_game_2013}.

\begin{example}
Consider the \IPG where player 1 solves $\max_{x^1}\{ 6x^1_1 + x^1_2 - 4x^1_1x^2_1 +6x^1_2x^2_2 : 3x^1_1 + 2x^1_2 \le 4, x^1\in\{0,1\}^2 \}$, and player 2 solves $\max_{x^2}\{ 4x^2_1 + 2x^2_2 -x^2_1x^1_1 -x^2_2x^1_2 : 3x^2_1 + 2x^2_2 \le 4, x^2\in\{0,1\}^2 \}$. The only \PNE is $(\bar{x}^1_1,\bar{x}^1_2)=(1,0)$, $(\bar{x}^2_1,\bar{x}^2_2)=(1,0)$ with a welfare of $S(\bar{x})=5$, $u^1(\bar{x}^1,\bar{x}^2)=2$, and $u^2(\bar{x}^2,\bar{x}^1)=3$.
\label{ex:KPGBoundDominance}
\end{example}
In the same fashion of \IEDS, we propose a family of inequalities that cuts off -- from each player's feasible set -- the strategies that never appear in a \PNE. Thus, from an \IPG instance $G$, we aim to derive an instance $G'$ where $\NESet^i$ replaces each player's feasible set $\mathcal{X}^i$. Note that, since all $\mathcal{X}^i$ are finite sets, all $\NESet^i$ are finite as well as the number of \PNEs.

\subsection{A Lifted Space}

Given the social welfare $S(x)$, we aim to find the \PNE maximizing it, namely, we aim to perform equilibria selection. In this context, the first urgent question is what space should we work in. Since mutually optimal strategies define \PNEs, a natural choice is to consider a space of all players' variables $x$.  As mentioned in the introduction, we assume the existence of a higher-dimensional (lifted) space where we linearize the non-linear terms in any $u^i(\cdot)$ via auxiliary variables $z$ and corresponding constraints (e.g., \citet{sherali_reformulation-linearization_1999,vielma_mixed_2015}).  Although our scheme holds for an arbitrary $f(x):\prod_{i=1}^n\mathcal{X}^i\rightarrow \mathbb{R}$ we can linearize to $f(x,z)$, we focus on $S(x)$ and the corresponding higher-dimensional $S(x,z)$ defined in the lifted space.  Let $\mathcal{L}$ be the set of (i.) linear constraints necessary to linearize the non-linear terms, and (ii.) the integrality requirements and bounds on the $z$ variables. The lifted space is then
\begin{align}
    \KSet= \{ (x^1,\dots,x^n,z) \in \mathcal{L}, x^i \in \mathcal{X}^i \text{ for any } i=1,\dots,n   \}.
    \label{eq:def:LiftedGame}
\end{align}
Any vector $x^1,\dots,x^n,z$ in (\ref{eq:def:LiftedGame}) corresponds to a unique strategy profile $x=(x^1,\dots,x^n)$, since $x$ induces $z$. $\KSet$ is then a set defined by linear constraints and integer requirements, and thus it is reasonable to deal with $\KSetConv$ and some of its projections. For brevity, let $\proj_{x}\KSetConv= \{x=(x^1,\dots,x^n) : \exists z \text{ s.t. } (x^1,\dots,x^n,z) \in \KSetConv \}$, and let $u^i(x^i,x^{-i})$ include the $z$ variables when working in the space of $\KSetConv$.

\subsection{Equilibrium Inequalities}
The integer points in $\proj_{x}(\KSetConv)$ encompass all the game's strategy profiles. However, we need to focus on $\EQFormulation =\{ (x^1,\dots,x^n,z) \in \KSetConv  : (x^1,\dots,x^n) \in \conv(\NESet) \}$, since projecting out $z$ yields the convex hull of \PNE profiles $\conv(\NESet)$. 
By definition $\EQFormulation $ is a polyhedron, and $\proj_{x^i}(\EQFormulation) = \conv(\NESet^i)$.
The role of $\EQFormulation$ is similar to the one of a perfect formulation for an integer program. As optimizing a linear function over a perfect formulation results in an integer optimum, optimizing a linear function $S(x,z)$ over $\EQFormulation$ results in a \PNE. For this reason, we call $\EQFormulation$ the \emph{perfect equilibrium formulation} for $G$. Also, the equivalent of the integrality gap in integer programming is the \POS, namely the ratio between the optimal value of $f(x,z)$ over $\KSetConv$ and $\EQFormulation$, respectively.
All considered, we establish the concept of \emph{equilibrium inequality}, a valid inequality for $\EQFormulation$.
\begin{definition}[Equilibrium Inequality]
Consider an \IPG instance $G$.
An inequality is an \emph{equilibrium inequality} for $G$ if it is a valid inequality for $\EQFormulation$.
\label{def:EqIneq}
\end{definition}

\paragraph{A Class of Equilibrium Inequalities. } We introduce a generic class of equilibrium inequalities that are linear in the space of $\KSetConv$. Consider any strategy $\tilde{x}^i \in \mathcal{X}^i$ for $i$: for any $i$'s opponents' strategy $x^{-i}$, $u^i(\tilde{x}^i,x^{-i})$ provides a lower bound on $i$'s payoff since $\tilde{x}^i\in \mathcal{X}^i$ (i.e., $\tilde{x}^i$ is a feasible point). Then,  $u^i(\tilde{x}^i,x^{-i}) \le u^i(x^i,x^{-i})$ holds for every player $i$. We introduce such inequalities in \cref{pro:Inequalities}.

\begin{proposition}
Consider an \IPG instance $G$. For any player $i$ and $\tilde{x}^i \in \mathcal{X}^i$, the inequality $u^i(\tilde{x}^i,x^{-i}) \le u^i(x^i,x^{-i})$ is an \emph{equilibrium inequality}.
\label{pro:Inequalities}
\end{proposition}
\iftoggle{ARXIV}{
    \begin{proof}
}{
    \proof{Proof of \cref{pro:Inequalities}.}
}
If a point $(\bar{x},\bar{z}) \in \EQFormulation$, then $\bar x \in \conv(\NESet)$. First, consider the case where $\bar{x} \in \ext(\conv(\NESet))$, namely $\bar{x} \in \NESet$ by definition. Assume $(\bar{x},\bar{z})$ violates the inequality associated with at least a player $i$, then, $u^i(\tilde{x}^i,\bar{x}^{-i}) > u^i(\bar{x}^i,\bar{x}^{-i})$. Therefore, $i$ can profitably deviate from $\bar{x}^i$ to $\tilde{x}^i$ under $\bar{x}^{-i}$, which contradicts $\bar{x} \in \NESet$ and $(\bar{x},\bar{z}) \in \EQFormulation$. Thus, no point $(\bar{x},\bar{z}) \in \EQFormulation$ with $\bar{x} \in \ext(\conv(\NESet))$ violates the inequality. Since we can represent any point $(\bar{x},\bar{z}) \in \EQFormulation$ as a convex combination of the extreme points of $\conv(\NESet)$, the proposition holds by iterating the previous reasoning for each extreme point in the support of $(\bar{x},\bar{z})$.
\iftoggle{ARXIV}{
    \end{proof}
}{
    \Halmos
    \endproof
}
A fundamental issue is whether the inequalities of \cref{pro:Inequalities} are sufficient to define the set $\EQFormulation$. By modulating the concept of closure introduced by  \citet{chvatal_edmonds_1973}, we prove this is indeed the case. We define the \emph{equilibrium closure} as the convex hull of the points in $\KSet$ satisfying the equilibrium inequalities of \cref{pro:Inequalities}.

\begin{theorem}
Consider an \IPG instance $G$ where $|\NESet| \neq 0$. Let the \emph{equilibrium closure} given by the equilibrium inequalities of \cref{pro:Inequalities} be 

$$
  P^{e}:= \conv \Big (\left\{ (x,z) \in \KSet : \begin{array}{l}
         u^i(\tilde{x}^i, x ^{-i}) \le u^i(x^i,x^{-i}) \\
        \forall \tilde{x}: \tilde{x}^i \in \mathcal{BR}(i,\tilde{x}^{-i}), \; i =1,\dots,n
  \end{array}\right\} \Big),
$$
where the equilibrium inequalities consider only the best-responses $\tilde{x}^i$ for any player $i$. Then, \begin{enumerate*} \item $P^{e}$ is a rational polyhedron, \item there exists no point $(x,z) \in \int(P^{e})$ such that  $x \in \mathbb{Z}^{nm}$, \item $P^{e} = \EQFormulation$ \end{enumerate*}.
\label{thm:Closure}
\end{theorem}
\iftoggle{ARXIV}{
    \begin{proof}
}{
    \proof{Proof of \cref{thm:Closure}.}
}
\emph{Proof of (i.) } The set $\KSet$ is finite since any $\mathcal{X}^i$ is finite, the number of best-responses and, correspondingly, of equilibrium inequalities, is finite.
Both equilibrium inequalities and the inequalities defining $\mathcal{X}^i$ have integer coefficients. Therefore, $P^{e}$ is a rational polyhedron.
\emph{Proof of (ii.) } Assume there exists a point $(\bar x, \bar z) \in \int(P^e)$ such that $\bar x \in \mathbb{Z}^{nm}$. By definition of Nash equilibrium, $\bar x \in \NESet$ since $(\bar x, \bar z)$ satisfies all the equilibrium inequalities in $P^e$. However, since $(\bar x, \bar z) \in \int(P^e)$, then no equilibrium inequality can be tight, contradicting the fact $\bar x$ is a \PNE. Therefore, there cannot exist any $(\bar x, \bar z) \in \int(P^e)$ such that $\bar x \in \mathbb{Z}^{nm}$. This also implies that all \PNEs lie on the boundary of $P^e$.
\emph{Proof of (iii.) } 
Since $P^e$ contains all the equilibrium inequalities generated by the players' best-responses, then any  $(\bar x, \bar z) \in \EQFormulation$ belongs to $P^e$ as of \cref{pro:Inequalities}, and $\EQFormulation \subseteq P^e$. 
Let $(\hat x, \hat z)$ be a point in $\ext(P^e)$. By definition, $(\hat x, \hat z)$ is an integer point, and it corresponds to a \PNE. Indeed, non-equilibria integer points cannot belong to $P^e$ since they would violate at least one equilibrium inequality associated with the players' best-responses.
Equivalently, for any $(\hat x, \hat z) \in \ext(P^e)$, its projection $\proj_x = \hat x$ is in $\NESet$.
Since all \PNEs are on the boundary of $P^e$, $P^{e} = \EQFormulation$ necessarily.
\iftoggle{ARXIV}{
    \end{proof}
}{
    \Halmos
    \endproof
}

Throughout the proof of \cref{thm:Closure}, we show that $P^e$ yields indeed the perfect equilibrium formulation $\EQFormulation$. Although the description of $P^e$ may contain an exponential number of possibly redundant equilibrium inequalities, it precisely describes the set of \PNEs in the lifted space. 
In \cref{ex:epsilon}, we showcase the construction $P^e$ via \cref{thm:Closure} for a small \IPG.

\begin{example}
\label{ex:epsilon}
Consider an \IPG where player 1 solves $\max_{x^1}\{ x^1_1 + 3x^1_2 + 7x^1_3 -6x^1_1x^2_1 +3x^1_2x^2_2 +2x^1_3x^2_3: 6x^1_1 + 4x^1_2 + 5x^1_3\le 7, x^1\in\{0,1\}^3 \}$, and player 2 solves $\max_{x^2}\{ 9x^2_1 + 9x^2_2 + 2x^2_3 -6x^2_1x^1_1 +5x^2_2x^1_2 +7x^2_3x^1_3 : 4x^2_1 + 2x^2_2 + 5x^2_3 \le 5, x^2\in\{0,1\}^3 \}$.
There are $4$ feasible strategies for each player $i$, namely, $(x^i_1,x^i_2,x^i_3)=(0,0,0)\lor(0,0,1)\lor(0,1,0)\lor(1,0,0)$. The $3$ \PNEs of this game are:
\begin{enumerate*}
\item $\bar{x}^1=(0,0,1)$ and $\bar{x}^2=(0,0,1)$ with $u^1(\bar{x}^1,\bar{x}^2)=9$ and $u^2(\bar{x}^2,\bar{x}^1)=9$, %
\item $\bar{x}^1=(0,0,1)$ and $\bar{x}^2=(0,1,0)$ with $u^1(\bar{x}^1,\bar{x}^2)=7$ and $u^2(\bar{x}^2,\bar{x}^1)=9$, %
\item $\bar{x}^1=(0,0,1)$ and $\bar{x}^2=(1,0,0)$ with $u^1(\bar{x}^1,\bar{x}^2)=7$ and $u^2(\bar{x}^2,\bar{x}^1)=9$.
\end{enumerate*}

We linearize the game by introducing $3$ variables $z_j \in \{0,1\}$ for any player's variable $j \in \{1,2,3\}$ such that $z_j=1$ if and only if $x^1_j=x^2_j=1$. We model these implications through the constraints $z_j \le x^i_j$ and $z_j \ge x^1_j + x^2_j -1$ for any player $i$ and variable $j$. Hence,
\begin{align}
\mathcal{K}= \left\{ 
		x^1 \in \{0,1\}^3, x^2\in \{0,1\}^3, z\in \{0,1\}^3  : 
		\begin{array}{l}
		    6x^1_1 + 4x^1_2 + 5x^1_3\le 7, \;
		    4x^2_1 + 2x^2_2 + 5x^2_3 \le 5 \\
		    z_j \le x^1_j, \; z_j \le x^2_j, \; z_j \ge x^1_j + x^2_j -1 \; \; \forall j \in \{1,2,3\}
		\end{array}
		\right\}. \nonumber
\end{align}
Correspondingly, the two players' utility functions in the linearized space are given by the two linear expressions $u^1(x^1,x^2)=x^1_1 + 3x^1_2 + 7x^1_3 -6z_1 +3z_2 +2z_3$ and $u^2(x^2,x^1)=9x^2_1 + 9x^2_2 + 2x^2_3 -6z_1 +5z_2 +7z_3$, respectively.

On the one hand, the best-response of player $1$ to any of player $2$'s feasible strategies is $\tilde{x}^1 = (0,0,1)$, i.e., $\mathcal{BR}(1,\tilde{x}^2) = \{ (0,0,1)\}$ for any feasible strategy $\tilde{x}^2$. The equilibrium inequality associated with $\tilde{x}^1=(0,0,1)$ is $7 + 2x^2_3 \le x^1_1 + 3x^1_2 + 7x^1_3 -6z_1 +3z_2 +2z_3$. The left-hand side of the inequality represents $u^1(\tilde{x}^1,x^{2})$, namely player 1's utility function evaluated on $\tilde{x}^1$.  
On the other hand, player $2$'s best-responses and the associated equilibrium inequalities are: \begin{enumerate*} 
\item $\tilde{x}^2=(1,0,0)$ with the inequality $9-6x^1_1 \le 9x^2_1 + 9x^2_2 + 2x^2_3 -6z_1 +5z_2 +7z_3$,
\item $\tilde{x}^2=(0,1,0)$ with the inequality $9+5x^1_2 \le9x^2_1 + 9x^2_2 + 2x^2_3 -6z_1 +5z_2 +7z_3$,
\item $\tilde{x}^2=(0,0,1)$ with the inequality $2+7x^1_3 \le 9x^2_1 + 9x^2_2 + 2x^2_3 -6z_1 +5z_2 +7z_3$.
\end{enumerate*}
Therefore, 
\begin{align}
  P^{e}= \conv \Big(\left\{ (x,z) \in \KSet \ : \begin{array}{l}
         7 + 2x^2_3 \le x^1_1 + 3x^1_2 + 7x^1_3 -6z_1 +3z_2 +2z_3\\
         9-6x^1_1 \le 9x^2_1 + 9x^2_2 + 2x^2_3 -6z_1 +5z_2 +7z_3 \\
         9+5x^1_2 \le9x^2_1 + 9x^2_2 + 2x^2_3 -6z_1 +5z_2 +7z_3 \\
         2+7x^1_3 \le 9x^2_1 + 9x^2_2 + 2x^2_3 -6z_1 +5z_2 +7z_3
  \end{array}\right\} \Big)  \nonumber.
\end{align}
By explicitly computing the above convex hull, we obtain
\begin{align}
P^e=  \left\{ (x,z) \ :\begin{array}{l}
         x^2_1 \ge 0,\; x^2_2 \ge 0,\; x^2_3\ge0,\; x^1_1 = 0,\; x^1_2 = 0,\; x^1_3 = 1, \\ x^2_1 + x^2_2 + x^2_3 = 1,\;
         z_1 = 0,\; z_2 = 0,\; x^2_1 + x^2_2 + z_3 = 1
  \end{array}\right\} \nonumber.
\end{align}
The projections onto the $x$ space of the extreme points of $P^e$ correspond to the $3$ \PNEs, and thus $P^e = \EQFormulation$. 
\end{example}

\section{The Cutting Plane Algorithm and its Oracle}
\label{sec:CompPNE}
If an oracle gives us $\EQFormulation$  in the form of a set of linear inequalities, then an optimal solution to $\max_{x^1,\dots,x^n,z}\{ f(x,z) : (x,z) \in \EQFormulation \}$ (i.e., a linear program) that is also an extreme point of $\EQFormulation$ is a \PNE for $G$ for any function $f(x,z)$. However, there are two major issues. First, $\EQFormulation \subseteq \KSetConv$, and $\KSetConv$ is a perfect formulation described by a possibly large number of inequalities. Second, retrieving $\EQFormulation$ through \cref{thm:Closure} may still require a large number of inequalities.
In practice, we actually do not need $\EQFormulation$ nor $\KSetConv$: a more reasonable goal is to get a polyhedron containing $\KSetConv$ over which we can optimize $f(x,z)$ efficiently and obtain an integer solution (i.e., $x \in \KSet$) that is also a \PNE. The first step is to obtain an integer solution. We could deploy branching schemes and known families of integer programming cutting planes, which are also equilibrium inequalities since they are valid for $\EQFormulation$. Equivalently, we can exploit a Mixed-Integer Programming (\MIP) solver to solve $\max_{x^1,\dots,x^n,z}\{ f(x,z) : (x,z) \in \KSet \}$. If the maximizer is a \PNE, the algorithm terminates. Otherwise, the second step is to cut off such maximizer, since it is not a \PNE, by separating at least an equilibrium inequality of \cref{pro:Inequalities}.

\paragraph{Equilibrium Separation Oracle. } Given a point $(\tilde x,\tilde z)$, for instance, the point returned by a \MIP solver, the central question is to decide if $\tilde x \in \NESet$,  and, if not, to derive an equilibrium inequality to cut off $(\tilde x,\tilde z)$. If we use the equilibrium inequalities of \cref{pro:Inequalities}, the process terminates in a finite number of iterations, since \cref{thm:Closure}. In the spirit of \citet{grotschel_ellipsoid_1981,karp_linear_1982}, we define a separation oracle for the equilibrium inequalities and $\EQFormulation$. The \emph{equilibrium separation oracle} solves the \emph{equilibrium separation problem} of \cref{def:SeparationProblem}.

\begin{definition}[Equilibrium Separation Problem]
Consider an \IPG instance $G$. Given a point $(\bar{x}, \bar{z})$, the \emph{equilibrium separation problem} is the task of determining that either: \begin{enumerate*} \item $(\bar{x}, \bar{z}) \in \EQFormulation$, or \item $(\bar{x}, \bar{z}) \notin \EQFormulation$ and return an equilibrium inequality violated by $(\bar{x}, \bar{z})$\end{enumerate*}. 
\label{def:SeparationProblem}
\end{definition}

\noindent \cref{Alg:ESO} presents our separation oracle for the inequalities of \cref{pro:Inequalities}. Given $(\bar{x}, \bar{z})$ and an empty set of linear inequalities $\phi$, the algorithm outputs either (i.) \texttt{yes} if $(\bar{x}, \bar{z}) \in \EQFormulation$, or (ii.) \texttt{no} and a set of violated equilibrium inequalities $\phi$ if $(\bar{x}, \bar{z}) \notin \EQFormulation$. The algorithm separates at most one inequality for any player $i$. By definition, $\bar{x}^i$ should be a best-response to be in a \PNE. Therefore, for any player $i$, the algorithm solves $\max_{x^i} \{ u^i(x^i,\bar{x}^{-i}) : A^i x^i \le b^i, x^i \in \mathbb{Z}^m \}$, where $\hat{x}^i$ is one of its maximizers. If $u^i(\bar{x}^i,\bar{x}^{-i}) = u^i(\hat{x}^i,\bar{x}^{-i})$, $\bar{x}^i$ is also a best-response.
However, if $u^i(\hat{x}^i,\bar{x}^{-i}) > u^i(\bar{x}^i,\bar{x}^{-i})$, the algorithm adds to $\phi$ an equilibrium inequality $u^i(\hat{x}^i,x^{-i})$ $\le u^i(x^i,x^{-i})$ violated by $(\bar{x}, \bar{z})$. After considering all players, if $|\phi|=0$, then $\bar{x}$ is a \PNE and the answer is \texttt{yes}. Otherwise, the algorithm returns \texttt{no} and $\phi \neq \emptyset$, i.e., at least an equilibrium inequality cutting off $(\bar{x}, \bar{z})$.

\begin{algorithm}[!ht]
	\caption{Equilibrium Separation Oracle \label{Alg:ESO}}
	\KwData{An \IPG instance $G$, a point $(\bar x, \bar z)$, and a set of cuts $\phi=\emptyset$.}
	\KwResult{ Either: (i.)
		\texttt{yes} if $(\bar{x}, \bar{z}) \in \EQFormulation$, or (ii.) \texttt{no} and $\phi$.}
	\For{$i\gets1$ \KwTo $n$}{	
    	$\hat{x}^i \gets \max_{x^i} \{ u^i(x^i,\bar{x}^{-i}) : A^i x^i \le b^i, x^i \in \mathbb{Z}^m \}$	\;
    	\If{$ u^i(\hat{x}^i,\bar{x}^{-i}) > u^i(\bar{x}^i,\bar{x}^{-i})$}
        	{
        	    	Add $u^i(\hat{x}^i,x^{-i})\le u^i(x^i,x^{-i})$ to $\phi$\;
        	}
	}
	\lIf{$|\phi|=0$} { \KwRet{\texttt{yes}} }
    \lElse { \KwRet{\texttt{no} and $\phi$} }
\end{algorithm}

\paragraph{ZERO Regrets. } We present our cutting plane algorithm \emph{ZERO Regrets} in \cref{Alg:ZERO}. 
The inputs are an instance $G$, and a  function $f(x)$, while the output is either the \PNE $\ddot x$ maximizing $f(x)$, or a certificate that no \PNE exists. Let $\Phi$ be a set of equilibrium inequalities, and $\mathcal{Q}=\max_{x^1,\dots,x^n,z} \{ f(x,z) : (x,z) \in \KSet, (x,z)  \; s.t. \; \Phi \}$. We assume $\mathcal{Q}$ is feasible and bounded. Otherwise, there would be no rationale behind getting a \PNE with an arbitrarily bad welfare. At each iteration, we compute an optimal solution $(\bar x, \bar z)$ of $\mathcal{Q}$. Afterwards, the equilibrium separation oracle (\cref{Alg:ESO}) evaluates $(\bar x, \bar z)$. If the oracle returns \texttt{yes}, then $\ddot{x}=\bar x$ is the \PNE maximizing $f(x)$ in $G$. Otherwise, the oracle returns a set $\phi$ of equilibrium inequalities cutting off $(\bar x, \bar z)$, and the algorithm adds $\phi$ to $\Phi$. Therefore, the process restarts by solving $\mathcal{Q}$ with the additional set of constraints. If at any iteration $\mathcal{Q}$ becomes infeasible, then $G$ has no \PNE. We remark that \cref{thm:Closure} implies both correctness and finite termination of \cref{Alg:ZERO}.
\begin{algorithm}[!ht]
	\caption{ZERO Regrets \label{Alg:ZERO}}
	\KwData{An \IPG instance $G$, and a function $f(x)$.}
	\KwResult{ Either: (i.) the \PNE $\ddot{x}$ maximizing $f(x)$, or (ii.) \texttt{no} \PNE }
	$\Phi =\{ 0 \le 1\}$, and $\mathcal{Q}=\max_{x^1,\dots,x^n,z} \{ f(x,z) : (x,z) \in \KSet, (x,z) \; s.t. \; \Phi \}$\;
	\While{true}{
	    \lIf{$\mathcal{Q}$ is infeasible}{\KwRet{\texttt{no} \PNE} } 
	    $(\bar x, \bar z)= \arg\max \mathcal{Q}$; $\phi=\emptyset$ \;
	    \If{\texttt{EquilibriumSeparationOracle$(G,(\bar x, \bar z), \phi)$} is \texttt{yes}}
	        {
	            \KwRet{$\ddot{x}=\bar x;$} \label{Alg:ZERO:ReturnPNE}
	        }
	    \lElse
	        {
	        add $\phi$ to $\Phi$
	        }
	}
\end{algorithm}
Although it is sufficient to add just one equilibrium inequality in $\phi$ cutting off the incumbent solution $(\bar x, \bar z)$, we expect that a good trade-off between $|\phi|=1$ and $|\phi|=n$ may speed up the convergence of the algorithm. This includes, for instance, also adding non-violated equilibrium inequalities. In \cref{ex:ZERORegrets}, we provide a toy example of \emph{ZERO Regrets}.
\begin{example}
Consider the game in \emph{Example 1} where player 1 solves $\max_{x^1}\{ 6x^1_1 + x^1_2 - 4x^1_1x^2_1 +6x^1_2x^2_2 : 3x^1_1 + 2x^1_2 \le 4, x^1\in\{0,1\}^2 \}$, and player 2 solves $\max_{x^2}\{ 4x^2_1 + 2x^2_2 -x^2_1x^1_1 -x^2_2x^1_2 : 3x^2_1 + 2x^2_2 \le 4, x^2\in\{0,1\}^2 \}$. 
As in \emph{Example 2}, to linearize the players' utility functions, we introduce two binary variables $z_1$ and $z_{2}$ equal to $1$ if both players select items $1$ and $2$, respectively. The linearization constraints are $z_1 \le x^1_1$, $z_1 \le x^2_1$, $z_1 \ge x^1_1 + x^2_1 - 1$, $z_2 \le x^1_2$, $z_2 \le x^2_2$, $z_2 \ge x^1_2 + x^2_2 - 1$. Thus, player 1's utility function is $6x^1_1 + x^1_2 - 4z_1 +6z_2$ and player 2's utility function is $4x^2_1 + 2x^2_2 - z_1 - z_2$. 
Correspondingly, problem $\mathcal{Q}$ maximizing the social welfare is
\begin{align}
\max_{(x^1, \;x^2, \;z)} \qquad & 6x^1_1 + x^1_2 + 4x^2_1 + 2x^2_2 -5z_1 + 5z_2  \nonumber \\
s.t. \qquad & 3x^1_1 + 2x^1_2 \le 4,\, 3x^2_1 + 2x^2_2 \le 4 \nonumber \\
& z_j \le x^1_j,\, z_j \le x^2_j,\, z_j \ge x^1_j + x^2_j - 1 \quad j=1,2. \nonumber \\
& x^1_j,\, x^2_j,\, z_j \in\{0,1\}  \quad j=1,2. \nonumber 
\end{align}
An optimal solution of the problem is $(\bar{x}^1_1,\bar{x}^1_2)=(1,0)$, $(\bar{x}^2_1,\bar{x}^2_2)=(0,1)$, $\bar{z}_1 = \bar{z}_2 = 0$. The social welfare is $8$, and the players' utility values are $6$ and $2$, respectively. However, this solution is not a \PNE. In fact, while a best-response to $\bar{x}^2$ for player 1 is $\bar{x}^1$, the best-response to $\bar{x}^1$ for player 2 is $(\hat{x}^2_1,\hat{x}^2_2)=(1,0)$ with an utility value of $3$.
Therefore, from player 2, we derive the equilibrium inequality
$4 - x^{1}_1  \le 4x^2_1 + 2x^2_2 -z_1 - z_2$
cutting off $(\bar{x}, \bar{z})$. 
By adding the equilibrium inequality to $\mathcal{Q}$, the optimal solution is then $(\bar{x}^1_1,\bar{x}^1_2)=(1,0)$, $(\bar{x}^2_1,\bar{x}^2_2)=(1,0)$, $\bar{z}_1 = 1$, $\bar{z}_2 = 0$, with utility values $2$ and $3$ and a welfare of $5$. Since $\bar{x}$ is a \PNE, the algorithm terminates by finding a \PNE with a \POS of $8/5$.
\label{ex:ZERORegrets}
\end{example}

\paragraph{Game-theoretical Interpretation.} We provide a straightforward game-theoretical interpretation of \emph{ZERO Regrets}. The algorithm acts as a central authority (e.g., a central planner) when optimizing $f(x,z)$ over $\KSet$, producing a solution that optimizes the welfare.
Afterward, it proposes the solution to each player, who evaluates it through the equilibrium separation oracle. The latter acts as a \emph{rationality} blackbox, in the sense that the oracle advises each player $i$ whether the proposed strategy is acceptable or not. In other words, the rationality blackbox tells the player $i$ if it should selfishly (and rationally) deviate to a better strategy, ignoring the central authority's advice. On the one hand, if the rationality blackbox says the solution is acceptable for player $i$,  then the player knows through the oracle that it should accept the proposed strategy. On the other hand, if at least one player $i$ refuses the proposed solution, the central authority should exclude such a solution and formulate a new proposal. Namely, it should cut off the non-equilibrium strategy and compute a new solution maximizing the welfare.

\subsection{Extensions}
\label{sec:extensions}
We showcase the flexibility of our algorithmic framework by proposing two extensions to \emph{ZERO Regrets}. Specifically, to address broader practical needs, we propose two extensions for enumerating \PNEs and computing approximate \PNEs.

\paragraph{Enumerating \PNEs. } We can easily tune \emph{ZERO Regrets} to enumerate all the \PNEs in $\NESet$ as follows. In \cref{Alg:ZERO:ReturnPNE} of \cref{Alg:ZERO}, instead of terminating and returning $\ddot x$, we memorize $\ddot x$ and add an (invalid) inequality cutting off the \PNE from $\EQFormulation$. Since all $x$ variables are integer-constrained, such inequality can be, for instance, an hamming distance from $\bar x$. The algorithm will possibly compute a new \PNE, cut it off (e.g., through a hamming distance constraint), and move the search towards the next equilibrium. Eventually, $\mathcal{Q}$ will become infeasible, thus certifying that the algorithm enumerated all the existing \PNEs. 

\paragraph{Approximating \PNEs. } An absolute $\epsilon$-\PNE is a \PNE where each player can deviate at most by a value $\epsilon$ for any best-response \citep{nisan_algorithmic_2008}, namely, where the regret for each player is at most $\epsilon$. Absolute $\epsilon$-\PNEs may be a reasonable compromise whenever no \PNE exists. Although any \PNE is an absolute $\epsilon$-\PNE with $\epsilon=0$, one may be interested in computing an absolute $\epsilon$-\PNE with an upper bound on $\epsilon$ while maximizing $f(x,z)$. We can adapt our algorithmic framework to compute an absolute $\epsilon$-\PNE as follows. We introduce a bounded continuous variable $\epsilon$ in $\mathcal{Q}$, and we let \cref{Alg:ESO} separate equilibrium inequalities in the form of $u^i(\hat{x}^i,x^{-i}) -  \epsilon \le u^i(x^i,x^{-i})$. Depending on the application of interest, one may still optimize the function $f(x,z)$ or  minimize $\epsilon$ without affecting the correctness of the algorithm. A similar modification enables the algorithm to handle relative $\epsilon$-\PNE, i.e., a profile of strategies where the payoff of each player's strategy is at least $\epsilon$ times the best-response payoff.  Given a constant $\epsilon$, the corresponding equilibrium inequalities are  $u^i(\hat{x}^i,x^{-i})  \epsilon \le u^i(x^i,x^{-i})$.

\section{Applications}
\label{sec:games}
We evaluate \emph{ZERO Regrets} on a wide range of problems from relevant works in the literature. We aim to provide a solid benchmark against the existing solution approaches and show the effectiveness of \emph{ZERO Regrets} in selecting and enumerating equilibria. The games we select stem from practical applications (e.g., competitive facility locations, network design) and methodological studies with the associated benchmark instances (e.g., games among quadratic programs).
Specifically, we consider the following games:
\begin{enumerate}
 \item The Knapsack Game (\KPG) \citep{Dragotto_2021_CNP,carvalho_computing_2017,carvalho_computation_2016}, where each player solves a binary knapsack problem. We select the equilibrium maximizing the social welfare, and we provide theoretical results on the complexity of deciding whether a \PNE exists. We also introduce two problem-specific equilibrium inequalities.
 \item The Network Formation Game (\NFG) \citep{chen_network_2006,anshelevich_price_2008} -- a paradigmatic game in \AGT with plenty of applications in network design --  where players seek to build a network through a cost-sharing mechanism. We select the equilibrium maximizing the social welfare.
 \item The Competitive Facility Location and Design game (\CFLD) \citep{cronert_equilibrium_2020,aboolian_competitive_2007}, where each player decides both the location of its facilities and their “design” (i.e., the facilities' features) while competing for customer demand. As in the \KPG and the \NFG, we focus on finding the \PNE maximizing the social welfare.
 \item The Quadratic \IPGs (\QIPGs) introduced by \citet{sagratella_computing_2016} and recently considered in \citet{schwarze_branch-and-prune_2022}, where each player optimizes a  (non-convex) quadratic function over box constraints and integrality requirements. As in the original papers, we focus on enumerating all the existing \PNEs, or determine that none exists.
\end{enumerate}
In what follows, we briefly describe the previous games and present the associated computational results\footnote{We performed our tests on an \emph{Intel Xeon Gold 6142} equipped with $128$GB of RAM and $8$ threads, employing $Gurobi$ 9.5 as \MIP solver for \cref{Alg:ZERO}.}. 

\subsection{The Knapsack Game}
The \KPG is an \IPG among $n$ players, where each player $i$ solves a binary knapsack problem with $m$ items in the form of
\begin{align}
    \max_{x^i} \Big\{ \sum_{j=1}^{m}p^{i}_jx^i_j + \sum_{k=1, k\neq i}^{n}\sum_{j=1}^{m} C^i_{k,j} x^i_j  x^k_j : \sum_{j=1}^{m} w^i_j x^{i}_j \leq b^i, x^i \in \{0,1\}^m \Big\}.
\label{eq:KPG}
\end{align}
As in the classical knapsack problem \citep{kellerer_knapsack_2004}, we assume that the profits $p^{i}_j$, weights $w^i_j$ and capacities  $b^i$ are in $\mathbb{Z}^+_0$. The selection of an item $j$ by a player $k \neq i$ impacts either negatively or positively the profit $p^{i}_j$ for player $i$ through integer coefficients $C^i_{k,j}$. Clearly, given the strategies of the other players $x^{-i}$, computing a corresponding best-response for player $i$ is \NPH. We can apply our algorithmic framework by linearizing the bilinear products $x^i_j x^k_j$ (for any $i,k,j$) with $\mathcal{O}(m n^2)$ auxiliary variables and additional constraints (see \cref{ex:epsilon}).
\citet{carvalho_computation_2016} introduced the game with $n=2$ and $p^i_j = 0$ $\forall j=1, \dots, m,\, i=1,2$. \citet{Dragotto_2021_CNP,carvalho_computing_2017} consider a more general game variant allowing $p^i_j$ and $w^i_j$ to take negative integer values. However, their algorithms focus on Mixed-Strategy equilibria and cannot perform exact equilibria selection.

In \cref{thm:IPGSigma2p}, we formally prove that deciding if a \KPG instance has a \PNE\xspace -- even with two players -- is \SigmaTwoPC in the polynomial hierarchy, matching the result of \citet{carvalho_computing_2017} for general \IPGs.
\begin{theorem}
 Deciding if a \KPG instance has a \PNE is a \SigmaTwoPC problem.
\label{thm:IPGSigma2p}
\end{theorem}
\noindent The proof, where we perform a reduction from the \SigmaTwoPC \emph{DeNegre Bilevel Knapsack Problem}  \citep{caprara_study_2014,denegre_thesis_2011, dellacroce_exact_2020}, is in the appendix.
Furthermore, we show that when at least one \PNE exists, the \POS and \POA can be arbitrarily bad.
\begin{proposition}
The \POA and the \POS in \KPG can be arbitrarily bad.
\label{prop:PoaPos}
\end{proposition}
\iftoggle{ARXIV}{
    \begin{proof}
}{
    \proof{Proof of \cref{prop:PoaPos}.}
}
Consider the following \KPG instance with $n=2$: player 1 solves the problem $\max_{x^1}\{ Mx^1_1 + x^1_2 - (M-2)x^1_1x^2_1 -x^1_2x^2_2  : 3x^1_1 + 2x^1_2 \leq 4, x^1\in\{0,1\}^2 \}$ where $M$ is an arbitrarily large value; player 2 solves $\max_{x^2}\{ 4x^2_1 + x^2_2 -x^2_1x^1_1 -x^2_2x^1_2 : 3x^2_1 + 2x^2_2 \leq 4, x^2\in\{0,1\}^2 \}$. The only \PNE is $(\bar{x}^1_1,\bar{x}^1_2,\bar{x}^2_1,\bar{x}^2_2)=(1,0,1,0)$, with $u^1(\bar{x}^1,\bar{x}^{2})=2$, $u^2(\bar{x}^2,\bar{x}^{1})=3$, $S(\bar{x})=5$. The maximum welfare $OSW=M+1$ is given by $(\hat{x}^1_1,\hat{x}^1_2,\hat{x}^2_1,\hat{x}^2_2)=(1,0,0,1)$, i.e., $OSW$ is arbitrarily large and there are no bounds on both the \POA and the \POS.
\iftoggle{ARXIV}{
    \end{proof}
}{
    \Halmos
    \endproof
}

\subsubsection{Strategic Inequalities. }
We further strengthen our cutting plane algorithm by introducing two classes of problem-specific equilibrium inequalities for the \KPG.

\paragraph{Strategic Dominance Inequalities.}  In the binary knapsack problem, a well-known hierarchy of dominance relationships exists among items, as we formalize in \cref{def:KPDominance}.
\begin{definition}[Dominance Rule]
Given two items $j$ and $j^\prime$ with profits $\bar{p}_j$ $\bar{p}_{j^\prime}$ and weights $w_j$, $w_{j^\prime}$, if $w_j \le w_{j^\prime}$ and $\bar{p}_j > \bar{p}_{j^\prime}$, then we say item $j$ \emph{dominates} item $j^\prime$.
\label{def:KPDominance}
\end{definition}

The above concept of dominance implies that, in any optimal knapsack solution, if one packs a dominated item $j^\prime$, then it should also pack item $j$. Otherwise, one would be able to improve the solution by selecting $j$ instead of $j^\prime$. This reasoning translates to the inequality $x_{j^\prime} \le x_j$, which is always valid for any optimal knapsack solution. We aim to extend this concept of dominance to the \KPG by incorporating the strategic interactions among players. In order to derive such inequalities, we reason about how, for any player $i$, the decisions of $x^{-i}$ affect the profits of $i$'s items. More formally, for any player $i$ and item $j$, let $p^{min}_j$, and $p^{max}_j$ be the minimum and maximum profit the strategies of the other players can induce, respectively. We claim the dominance rule of \cref{def:KPDominance} extends to the one of \cref{prop:StrictDom} in the \KPG.

\begin{proposition}
For each player $i$, if the dominance rule applies for two items $j$ and $j^\prime$ with $\bar{p}_j = p^{min}_j$ and $\bar{p}_{j^\prime} = p^{max}_{j^\prime}$, then the inequality $x^i_{j^\prime} \le x^i_{j}$ is an \emph{equilibrium inequality}.
\label{prop:StrictDom}
\end{proposition}
\iftoggle{ARXIV}{
    \begin{proof}
}{
    \proof{Proof of \cref{prop:StrictDom}.}
}
Since all best-responses of player $i$ cannot select the dominated item $j^\prime$ without selecting item $j$ for any $x^{-i}$, the claim holds.
\iftoggle{ARXIV}{
    \end{proof}
}{
    \Halmos
    \endproof
}

We denote the inequalities of \cref{prop:StrictDom} as \emph{Strategic Dominance Inequalities}. We further extend the previous reasoning to derive other forms of dominance inequalities by evaluating how the strategic interaction (i.e., the items that the other players select) affects the items' profits for each player $i$. In other words, we derive equilibrium inequalities that incorporate the strategic interaction by including the variables of multiple players. For instance, consider the case with two players. If the profits of two items $j$ and $j^\prime$ for player $1$ fulfill the dominance rule when player $2$ selects item $j$ and does not select item $j^\prime$, then 
$$ x^1_{j^\prime} \le  x^1_{j} + (1 - x^2_j) + x^2_{j^\prime}$$
is an equilibrium inequality. Namely, if there exists a \PNE with $x^2_j = 1$ and $x^2_{j^\prime} = 0$, the dominance rule between item $j$ and $j^\prime$ applies for player $1$, otherwise the inequality is not binding.

\paragraph{Strategic Payoff Inequalities. } We introduce a second class of \emph{strategic} inequalities by exploiting two observations on the knapsack problem. For any player $i$, the strategy of all-zeros $x^i=(0,\dots,0)$ is always feasible under the packing constraint, and yields a payoff of $0$. Therefore, for any player $i$ and item $j$, if $p^i_j + \sum_{k=1, k \neq i}^n C^i_{k,j} <0$, player $i$ may not select $j$ depending on its opponent choices $x^{-i}$. More generally, let $\mathcal{S}^{i}_j$ be the \emph{interaction set} of $i$'s opponents inducing a negative profit for item $j$, namely, a set so that
\begin{align}
p^i_j + \sum_{k \in \mathcal{S}^{i}_j} C^i_{k,j} <0.
\label{eq:StrategicPayoff1}
\end{align}
The interaction is \emph{minimal} if, for any proper subset $\bar{\mathcal{S}}^{i}_j$ of $\mathcal{S}^{i}_j$, then $p^i_j + \sum_{k \in \bar{\mathcal{S}}^{i}_j} C^i_{k,j} >0$. The inequality (\ref{eq:StrategicPayoff1}) implies that if $x^k_j=1$ for any $k \in \mathcal{S}^{i}_j$, then $x^i_j=0$. In general, this implies that for any interaction set, the inequality
$$
x^i_j + \sum_{k \in \mathcal{S}^{i}_j} x^k_j \le |\mathcal{S}^{i}_j|
$$
is an equilibrium inequality. We define the latter inequality as \emph{Strategic Payoff Inequality}. In practice, the inequalities generated by minimal interaction sets are stronger than those generated by non-minimal interaction sets since they generally involve more variables. 
Clearly, the effort to separate and include all the previous strategic inequalities may not be negligible when $n$ and $m$ increase. In practice, at each iteration of \cref{Alg:ZERO}, we separate and add to $\mathcal{Q}$ only the inequalities violated by the incumbent solution $(\bar x, \bar z)$.

\paragraph{Computational Results. } We generate \KPG instances with $n= 2, 3$ and $m= 25, 50, 75,$ $100$, and with $p^{i}_j$ and $w^i_j$ being random integers uniformly distributed in $[1, 100]$ for any $i$. The knapsack capacities $b^i$ are equal to $0.2\sum_{j=1}^{m} w^i_j$, $0.5\sum_{j=1}^{m} w^i_j$, or $0.8\sum_{j=1}^{m} w^i_j$, respectively. For what concerns the strategic interaction, we focus on three different distributions for the integer interaction coefficients $C^i_{k,j}$. For any player $i$, they can be: \emph{A)} equal and uniformly distributed in $[1,100]$, or \emph{B)} random and uniformly distributed in $[1,100]$, or \emph{C)} random and uniformly distributed in $[-100,100]$. In \cref{tab:KPG}, we present the results for the $72$ resulting instances. For any given number of players $n$, items $m$ and distribution of coefficients $C^i_{k,j}$ ($(n, m, d$)), we report the performance over the $3$ instances with different capacities, in terms of:
\begin{enumerate*}
\item the average number of Equilibrium Inequalities of \cref{pro:Inequalities} added (\emph{\#EI}), 
\item the average number of Strategic Payoff Inequalities (\emph{\#EI\_P}) (which we only compute for the instances with distribution \emph{C}), 
\item the average number of Strategic Dominance Inequalities (\emph{\#EI\_D}), 
\item the average computational time (\emph{Time}),
\item the average computational time (\emph{Time-\nth{1}}) to find the first \PNE (if any),
\item the average \POS (\POS) for the best \PNE (if any), and
\item the number of time limit hits (\emph{Tl})
\end{enumerate*}. The average values in \emph{\#EI}, \emph{\#EI\_P}, \emph{\#EI\_D}, \emph{Time} and \emph{Time-\nth{1}} also consider the instances where we hit the time limit\footnote{ We remark that the same observation holds on all our experiments.}, which we set to $1800$ seconds.
\emph{ZERO Regrets} solves almost all instances with $n=2$, regardless of the type of strategic interaction. Both running times and the number of equilibrium inequalities are significantly modest for a \SigmaTwoP game. The \POS is generally low and increases with distribution $C$ due to the nature of the complex strategic interactions stemming from both negative and positive $C^i_{k,j}$ coefficients. We remind that a \POS close to $1$ does not mean the instance is computationally “easy”. On the contrary, a $\POS \approx 1$ highlights the existence of a high-quality \PNE (i.e., with a welfare close to the one of the $OSW$) and also provides further evidence concerning the urgency of selecting efficient \PNEs.
\emph{ZERO Regrets} performs robustly even in large instances, establishing a significant computational advantage over the previously developed approaches in the literature. \citet{carvalho_computing_2017,Dragotto_2021_CNP} consider up to $m=40$ items with $n=3$ by just computing an equilibrium, while \citet{cronert_equilibrium_2020} only perform equilibria selection with $m<5$. 

\begin{table}[!ht]
\setlength\tabcolsep{1em}
\centering
\resizebox{0.9\textwidth}{!}{%
\begin{tabular}{>{\bf}c@{\hspace{3em}}lllllll} 
\toprule
\textbf{($\bm n$, $\bm m$, $\bm d$)}    & \textbf{\#EI}  & \textbf{\#EI\_P} & \textbf{\#EI\_D} & \textbf{Time} & \textbf{Time-\nth{1}}    & \textbf{\POS}  & \textbf{Tl}   \\
\midrule
(2, 25, A)  & 14.67  & 0.00     & 3.00      & 0.06    & 0.05   & 1.07 & 0/3 \\
(2, 25, B)  & 17.33  & 0.00     & 3.67      & 0.12    & 0.09   & 1.02 & 0/3 \\
(2, 25, C)  & 29.33  & 9.67     & 7.67      & 0.39    & 0.04   & 1.06 & 0/3 \\
(2, 50, A)  & 20.00  & 0.00     & 2.67      & 0.21    & 0.21   & 1.02 & 0/3 \\
(2, 50, B)  & 26.67  & 0.00     & 19.67     & 0.51    & 0.39   & 1.01 & 0/3 \\
(2, 50, C)  & 72.67  & 27.00    & 11.33     & 6.34    & 0.92   & 1.08 & 0/3 \\
(2, 75, A)  & 38.00  & 0.00     & 31.00     & 0.60    & 0.44   & 1.00 & 0/3 \\
(2, 75, B)  & 100.67 & 0.00     & 34.00     & 8.35    & 5.71   & 1.02 & 0/3 \\
(2, 75, C)  & 112.67 & 38.33    & 67.00     & 47.75   & 3.96   & 1.08 & 0/3 \\
(2, 100, A) & 25.33  & 0.00     & 14.67     & 0.76    & 0.58   & 1.01 & 0/3 \\
(2, 100, B) & 205.33 & 0.00     & 79.67     & 220.42  & 143.45 & 1.01 & 0/3 \\
(2, 100, C) & 697.33 & 55.33    & 119.67    & 1205.29 & 11.33  & 1.05 & 2/3 \\
(3, 25, A)  & 31.00  & 0.00     & 9.33      & 0.21    & 0.21   & 1.01 & 0/3 \\
(3, 25, B)  & 44.00  & 0.00     & 14.67     & 0.33    & 0.33   & 1.02 & 0/3 \\
(3, 25, C)  & 91.00  & 29.67    & 33.67     & 29.78   & 5.64   & 1.26 & 0/3 \\
(3, 50, A)  & 95.00  & 0.00     & 24.33     & 18.39   & 11.68  & 1.03 & 0/3 \\
(3, 50, B)  & 206.00 & 0.00     & 44.33     & 626.45  & 167.01 & 1.01 & 1/3 \\
(3, 50, C)  & 148.00 & 63.00    & 224.67    & 382.24  & -      & -    & 0/3 \\
(3, 75, A)  & 64.00  & 0.00     & 119.00    & 4.65    & 2.07   & 1.02 & 0/3 \\
(3, 75, B)  & 278.00 & 0.00     & 92.67     & 982.97  & 272.69 & 1.01 & 1/3 \\
(3, 75, C)  & 173.00 & 87.33    & 319.67    & 658.77  & -      & -    & 1/3 \\
(3, 100, A) & 261.00 & 0.00     & 144.67    & 1200.65 & 666.13 & 1.00 & 2/3 \\
(3, 100, B) & 479.00 & 0.00     & 168.33    & tl & -      & -    & 3/3 \\
(3, 100, C) & 184.00 & 171.00   & 1019.67   & 1200.31 & -      & -    & 2/3 \\
\bottomrule \vspace{-\aboverulesep}
\end{tabular} %
}
\caption{Results overview for the \KPG. The complete tables of results are in the Appendix (\cref{tab:KPG_Full2} and \cref{tab:KPG_Full3}).}
\label{tab:KPG}
\end{table}

\subsection{The Network Formation Game}
Network design games are paradigmatic problems in Algorithmic Game Theory \citep{chen_network_2006,anshelevich_price_2008,nisan_algorithmic_2008}. 
Their natural application domain is often the one of computer networks and the Internet itself, where several selfish agents opportunistically decide how to share a scarce resource, for instance, the bandwidth. \citet{tardos_network_2004} accurately claims that the impact and future of the complex technology we develop through the Internet critically depend on the ability to balance the diverse needs of the selfish agents in the network. We consider a (weighted) \NFG -- similar to the one of \citet{chen_network_2006} -- where $n$ players are interested in building a computer network. Let $G(V, E)$ be a directed graph representing a network layout, where $V$, $E$ are the sets of vertices and edges, respectively. Each edge $(h,l) \in E$ has a construction cost $c_{hl} \in \mathbb{Z}^+$, and each player $i$ wants to connect an origin $s^i$ with a destination $t^i$ while minimizing its construction costs. A cost-sharing mechanism determines the cost of each edge $(h,l)$ for a player as a function of the number of players crossing $(h,l)$. Arguably, the most common and widely-adopted mechanism is the \emph{Shapley cost-sharing mechanism}, where players using $(h,l)$ equally share its cost $c_{hl}$. The goal is to find a \PNE minimizing the sum of construction costs for each player or determine that no \PNE exists. 
We model the \NFG as an \IPG as follows. For any player $i$ and edge $(h,l)$, let the binary variables $x^i_{hl}$ be $1$ if $i$ uses the edge. We employ standard flow constraints to model a path between $s^i$ and $t^i$.
For conciseness, we represent these constraints and binary requirements with a set $\mathcal{F}^i$ for each $i$. Thus, each player $i$ solves:
\begin{align}
    \min_{x^i} \{ \sum_{(h,l)\in E} \frac{c_{hl}x^i_{hl}}{\sum_{k=1}^{n} x^k_{hl}}: x^i \in \mathcal{F}^i \}.
   \label{eq:NFG}
\end{align}
For any player $i$, the cost contribution of each edge $(h,l)$ to the objective is not linear in $x$ and may not be defined for some choices of $x$ (i.e.,  $\sum_{k=1}^{n} x^k_{hl}=0$). However, we can linearize the fractional terms and eliminate the indefiniteness. For instance, consider a game with $n=3$ and the objective of player $1$. Let the binary variable $z^{j, \dots, k}_{hl}$ be $1$ if only players $j, \dots, k$ select the edge $(h,l)$. Then, $x^1_{hl} = z^1_{hl} + z^{12}_{hl} + z^{13}_{hl} + z^{123}_{hl}$, $x^2_{hl} = z^2_{hl} + z^{12}_{hl} + z^{23}_{hl} + z^{123}_{hl}$, $x^3_{hl} = z^3_{hl} + z^{13}_{hl} + z^{23}_{hl} + z^{123}_{hl}$ along with a clique constraint $ z^1_{hl} + z^2_{hl} + z^3_{hl} + z^{12}_{hl} + z^{13}_{hl} + z^{23}_{hl} + z^{123}_{hl} \le 1$. The term for edge $(h,l)$ in the objective of player $1$ is then $c_{hl}z^1_{hl} + \frac{c_{hl}}{2}(z^{12}_{hl} + z^{13}_{hl}) + \frac{c_{hl}}{3}z^{123}_{hl}$. In our tests, we consider the general weighted \NFG \citep{chen_network_2006}, where each player $i$ has a weight $w^i$, and the cost share of each selected $(h,l)$ is $w^i c_{hl}$ divided by the weights of all players using $(h,l)$. Specifically, we consider the $3$-player weighted \NFG, where a \PNE may not exist, and selecting one if multiple equilibria exist is generally an \NPH problem \citep{anshelevich_price_2008,chen_network_2006}.

\paragraph{Computational Results. } In order to tackle challenging instances, we consider the \NFG with $n = 3$ on grid-based directed graphs $G(V, E)$, where each $i$ has to cross the grid from left to right to reach its destination. Compared to a standard grid graph, we randomly add some edges between adjacent layers to increase the number of paths, and to facilitate the interaction among players. The instances are so that $|V| \in [50,500]$, and the costs $c_{hl}$ for each edge $ (h,l)$ are random integers uniformly distributed in $[20, 100]$. 
We consider three distributions of player's weights: \begin{enumerate*}
\item the Shapely-mechanism with $w^1=w^2=w^3=1$, where a \PNE always exists, yet selecting the most efficient \PNE is \NPH, or
\item $w^1=0.6$, $w^2=0.2$, and $w^3=0.2$, or
\item $w^1=0.45$, $w^2=0.45$, and $w^3=0.1$.
\end{enumerate*}
\cref{tab:NFG} reports the results, where we average over the distributions of the players' weights. For each graph, the table reports the graph size ($|V|,|E|$), whereas the other columns have the same meaning of the ones of \cref{tab:KPG}. Our algorithm effectively solves all the instances but $3$ within a time limit of $1800$ seconds and consistently selects high-efficiency \PNEs. Further, our algorithm finds the first \PNE in considerably limited computing times. Generally, the previous literature does not consider this problem from a computational perspective, but only provides theoretical and possibly pessimistic bounds on the \POS and \POA. Nevertheless, we can compute efficient \PNEs even in large-size graphs (i.e., $\POS \approx 1$), with a limited number of equilibrium inequalities and modest running times, showing the practical effectiveness of our algorithm in a paradigmatic \AGT problem.

\begin{table}[!ht]
\centering
\resizebox{\textwidth}{!}{%
\begin{tabular}{>{\bf}c@{\hspace{1em}}r@{\hspace{1em}}r@{\hspace{1em}}r@{\hspace{1em}}r@{\hspace{3em}}r@{\hspace{1em}}>{\bf}c@{\hspace{1em}}r@{\hspace{1em}}r@{\hspace{1em}}r@{\hspace{1em}}r@{\hspace{1em}}r@{\hspace{1em}}}
\toprule
\textbf{($\bm{|V|}$, $\bm{|E|}$)} & \textbf{\#EI} & \textbf{Time} & \textbf{Time-\nth{1}}& \textbf{\POS} & \textbf{Tl} & \textbf{($\bm{|V|}$, $\bm{|E|}$))} & \textbf{\#EI} & \textbf{Time} &  \textbf{Time-\nth{1}}& \textbf{\POS} & \textbf{Tl} \\
\midrule
(50, 99)   & 6.00  & 0.04  & 0.04  & 1.12 & 0/3  & (300, 626)  & 21.00  & 12.11   & 2.64   & 1.04 & 0/3  \\
(100, 206) & 2.33  & 0.05  & 0.04  & 1.00 & 0/3  & (350, 730)  & 19.00  & 13.92   & 7.42   & 1.01 & 0/3  \\
(150, 308) & 6.00  & 0.64  & 0.25  & 1.01 & 0/3  & (400, 822)  & 248.67 & 694.95  & 228.69 & 1.08 & 1/3  \\
(200, 416) & 11.67 & 3.28  & 1.11  & 1.06 & 0/3  & (450, 934)  & 394.67 & 1199.98 & 2.61   & 1.11 & 2/3  \\
(250, 517) & 64.67 & 63.50 & 16.07 & 1.02 & 0/3  & (500, 1060) & 35.67  & 87.07   & 7.25   & 1.00 & 0/3  \\
\bottomrule \vspace{-\aboverulesep}  \\
\end{tabular} %
}
\caption{Results overview for the \NFG. The complete table of results is in the Appendix (\cref{tab:NFG_Full}).}
\label{tab:NFG}
\end{table}

\subsection{The Competitive Facility Location and Design Game}
The \CFLD \citep{aboolian_competitive_2007} is a game where sellers (players) compete for the demand of customers located in a given geographical area. Each seller makes two fundamental choices: where to open its selling facilities and the product assortment of such facilities, i.e., their “design”. Symmetrically, the customers select their favorite facilities depending on the relative distance from a facility and its attractiveness in terms of design.
We consider a variant recently presented by \citet{cronert_equilibrium_2020}, where $n$ competitors simultaneously choose the location and design of their facilities. Let $L$ be the set of potential facility locations, $J$ be the set of customers, and let $R_l$ denote the set of design alternatives for each location $l \in L$. Each player $i$ has an available budget $B^i$ and incurs in a fixed cost $f_{lr}^i$ when opening a facility at location $l \in L$ and with the design $r \in R_l$. Each player $i$ acquires a share of the demand $w_j$ of a customer $j \in J$ according to a utility $u_{ljr}^i$, whose value depends upon the distance of customer $j$ from a facility in location $l$ and the design choice of such facility (see \citet{cronert_equilibrium_2020} for more details). The \CFLD formulates as an \IPG where each player $i$ solves
\begin{maxi!}
	{x^i}{ 
	\sum_{j\in J} w_j \frac{\sum_{l \in L} \sum_{r \in R_l} u_{ljr}^ix^i_{lr}}
     {\sum_{k=1}^n \sum_{l \in L} \sum_{r \in R_l} u_{ljr}^kx^k_{lr}}
	\protect\label{eq:CFLD:Obj}}
	{\label{eq:CFLD}}{}
	\addConstraint{\sum_{l \in L} \sum_{r \in R_l} f_{lr}^ix^i_{lr}}{ \le B^i \protect\label{eq:CFLD:Budget}, \quad}{}{}
	\addConstraint{\sum_{r \in R_l} x^i_{lr} }{\le 1 \protect\label{eq:CFLD:sos} \quad \quad }{\forall l \in L,}
	\addConstraint{ x^i_{lr} \in \{0, 1\}}{ \protect\label{eq:CFLD:Vars} \quad }{\forall l \in L, \forall r \in R_l.}
\end{maxi!}
  
\noindent The binary variable $x^i_{lr}$ is $1$ if player $i$ opens a facility in the location $l \in L$ with a design $r \in R_l$.
The objective function \eqref{eq:CFLD:Obj} represents the share of customer demands player $i$ maximizes, the constraint \eqref{eq:CFLD:Budget} is the budget constraint for player $i$, and the constraints \eqref{eq:CFLD:sos} enforce that player $i$ can open only one facility in a location $l$. As in the \NFG, the objective is not linear in $x$, and the denominator can be zero; however, we can linearize it with tailored fractional programming techniques.

\paragraph{Computational Results. }
We test \emph{ZERO Regrets} on a representative set of instances from \citet{cronert_equilibrium_2020} to which we refer for the details concerning the distributions of locations and customers, and the entries $w_j$, $u_{ljr}$, $f_{lr}$.  The resulting $64$ instances with $n=2,3$ have $50$ locations and $50$ customers, with budgets $B^1 \in [10, 40]$ and  $B^2 = B^1, B^1 +10, \dots, 100$, $B^3 = 10$. \cref{tab:CFLD} summarizes the results, where we aggregate and average over the values of $B^1$. We benchmark our results against the performance of \emph{eSGM-WM} from \citet[Table 2]{cronert_equilibrium_2020}, which ran on a machine with similar hardware characteristics. Although the authors compute both pure and mixed welfare-maximizing equilibria, we focus on computing only the welfare-maximizing \PNE. Generally, \emph{ZERO Regrets} outperforms algorithm \emph{eSGM-WM} even in instances where only \PNEs exist. Our algorithm solves $62$ instances over $64$ within a time limit of $3600$ seconds. The running times of \emph{ZERO Regrets} are sensibly limited compared to those of \emph{eSGM-WM}, and never hit the time limit on the instances with $n=2$. Occasionally, the running times are dramatically smaller, e.g., the instance with $n=2$, $B^1 = 40$, $B^2 = 80$ where only one \PNE exists: our algorithm finds the most efficient \PNE in about $1636$ seconds, while \emph{eSGM-WM} requires $163315$ seconds.

\begin{table}[!ht]
\centering
\resizebox{\textwidth}{!}{%
\begin{tabular}{>{\bf}c@{\hspace{1em}}r@{\hspace{1em}}r@{\hspace{1em}}r@{\hspace{1em}}r@{\hspace{3em}}r@{\hspace{1em}}>{\bf}c@{\hspace{1em}}r@{\hspace{1em}}r@{\hspace{1em}}r@{\hspace{1em}}r@{\hspace{1em}}r@{\hspace{1em}}}
\toprule
\textbf{($\bm{n}$, $\bm{B^1}$)} & 
\textbf{\#EI} & 
\textbf{Time} & 
\textbf{Time-\nth{1}}& 
\textbf{\POS} & 
\textbf{Tl} & 
\textbf{($\bm{n}$, $\bm{B^1}$)} & 
\textbf{\#EI} & 
\textbf{Time} &  
\textbf{Time-\nth{1}}& 
\textbf{\POS} & 
\textbf{Tl} \\
\midrule
(2, 10) & 4.00  & 1.76   & 1.76   & 1.01 & 0/10 & (3, 10) & 6.30  & 3.56    & 2.11   & 1.02 & 0/10 \\
(2, 20) & 5.11  & 7.39   & 3.03   & 1.01 & 0/9  & (3, 20) & 9.00  & 27.72   & 7.29   & 1.03 & 0/9  \\
(2, 30) & 9.25  & 339.35 & 59.39  & 1.03 & 0/8  & (3, 30) & 18.38 & 754.84  & 555.78 & 1.05 & 1/8  \\
(2, 40) & 16.40 & 682.00 & 294.92 & 1.07 & 0/5  & (3, 40) & 25.20 & 1863.92 & 739.63 & 1.06 & 1/5  \\
\bottomrule \vspace{-\aboverulesep}
\end{tabular} %
}
\caption{Results overview for the \CFLD from the instances of \citet[Table 2]{cronert_equilibrium_2020}. The complete table of results is in the Appendix (\cref{tab:CFLD_Full}).}
\label{tab:CFLD}
\end{table}

\subsection{The Quadratic Game}
The \QIPG is a simultaneous non-cooperative \IPG introduced by \citet{sagratella_computing_2016}, where each player $i$ solves the problem
\begin{align}
\min_{x^i} \{ \frac{1}{2}(x^i)^\top Q^i x^i + (C^i x^{-i})^\top x^i + (c^i)^\top x^i : LB \le x^i \le UB, \; x^i \in \mathbb{Z}^m \}.
\label{eq:QuadratiGame}
\end{align}
Specifically, each player $i$ controls $m$ integer variables bounded by the vectors $LB$ and $UB$. The strategic interaction involves the term $(C^i x^{-i})^\top x^i$, while the linear and quadratic terms solely depend on each player's choices. While \citet{sagratella_computing_2016} considers only instances with positive-definite $Q^i$ matrices (i.e., the problem is convex in $x^i$ for any $i$), \citet{schwarze_branch-and-prune_2022} consider arbitrary matrices $Q^i$ (i.e., non-convex objectives). In particular, the latter generalizes the former by dropping the convexity requirement w.r.t. $x^i$ on the payoffs $u^i(x^i,x^{-i})$. In contrast with the aforementioned applications, we let the \MIP solver manage the linearization of the quadratic terms in each player's payoff in order to fully integrate \emph{ZERO Regrets} with the features of the existing \MIP technology. As in \citet{sagratella_computing_2016,schwarze_branch-and-prune_2022}, we setup \emph{ZERO Regrets} to enumerate all \PNEs or to certify that no \PNE exists.

\paragraph{Computational Results.} We test our algorithm on both convex and non-convex benchmarks of the \QIPG. 
First, we consider the \QIPG from \citet{schwarze_branch-and-prune_2022}, and test our algorithm on the same instance set. We refer to the original paper for the details on instance generation. Besides the bounds on the $x^i$ variables, these instances also include $m$ non-redundant linear inequalities $A^i x^i \le b^i$ for each player $i$. \cref{tab:Quad_SS} reports an overview of the results with a similar notation to the one of the previous tables.  In the first column, we report the tuple $(n, m, t)$, where $t$ is either $C$ when each player's problem is convex or $NC$ otherwise. We additionally report the average number of \PNEs in the column $\#EQs$. We solve all the $56$ instances in less than $416$ seconds globally, with the most computationally-difficult instance requiring $56$ seconds.
Similar to the previous tests, our algorithm strongly outperforms the one of \citet{schwarze_branch-and-prune_2022}.  
Specifically, their algorithm runs out of time in $25$ instances (time limit of $3600$ seconds) and solves the remaining $31$ instances with non-negligible computational times (i.e., about $1302$ seconds on average).

\begin{table}[!ht]
\setlength\tabcolsep{0.7em}
\centering
\resizebox{0.67\textwidth}{!}{%
\begin{tabular}{>{\bf}crrrrrrr}
\toprule
\textbf{($\bm{n}, \bm{m}, \bm{t}$)} & 
\textbf{\#EI} & 
\textbf{\#EQs} &
\textbf{Time} & 
\textbf{Time-\nth{1}} &
\textbf{\POS} &
\textbf{\POA} &
\textbf{Tl} \\
\midrule
(2, 2, C)  & 14.00  & 1.67  & 0.35  & 0.18 & 1.17 & 1.43 & 0/4 \\
(2, 3, C)  & 26.75  & 1.60  & 1.13  & 0.42 & 5.56 & 6.10 & 0/8 \\
(2, 4, C)  & 35.00  & 1.00  & 3.40  & 1.19 & 1.32 & 1.32 & 0/4 \\
(2, 5, C)  & 56.00  & 1.50  & 35.88 & 7.72 & 2.19 & 4.74 & 0/4 \\
(3, 2, C)  & 42.00  & 2.00  & 1.21  & 0.65 & 1.83 & 2.40 & 0/4 \\
(3, 3, C)  & 108.75 & 1.00  & 35.62 & 7.09 & 3.79 & 3.79 & 0/4 \\
(2, 2, NC) & 16.00  & 1.33  & 0.43  & 0.30 & 2.20 & 2.20 & 0/4 \\
(2, 3, NC) & 20.25  & 1.75  & 1.07  & 0.40 & 1.33 & 1.68 & 0/8 \\
(2, 4, NC) & 19.50  & 1.00  & 1.51  & 0.82 & 1.02 & 1.02 & 0/4 \\
(2, 5, NC) & 28.50  & 1.67  & 14.03 & 2.38 & 1.31 & 1.53 & 0/4 \\
(3, 2, NC) & 30.00  & 2.67  & 1.03  & 0.37 & 1.03 & 1.60 & 0/4 \\
(3, 3, NC) & 46.50  & 1.33  & 6.11  & 2.43 & 1.44 & 1.45 & 0/4 \\
\bottomrule \vspace{-\aboverulesep} 
\end{tabular} %
}
\caption{Results overview for the \QIPG from the instances of \citet{schwarze_branch-and-prune_2022}. The complete table of results is in the Appendix (\cref{tab:Quad_SS_Full}).}
\label{tab:Quad_SS}
\end{table}

To get a broader perspective, we also consider the convex instances of the game generated according to the scheme proposed in \citet{sagratella_computing_2016}. The matrices $Q^i$ and $C^i$ are a random positive-definite and a random matrix, respectively, with rational entries in the range $[-25, 25]$, while the entries of $c^i$ are integer numbers in the range $[-5, 5]$. We generate our instances with $n \in [1,6]$, $m \in [1,10]$, $LB \in [-1000,0]$, and $UB \in[5,1000]$  similarly to \citet{sagratella_computing_2016}.
We report the average results in \cref{tab:Quad_S}, where we aggregate for $n$. 
\emph{ZERO Regrets} finds the first \PNE in less than a second on average, and manages to solve any instance in less than $12$ seconds, even when more than $10$ \PNEs exist (see \cref{tab:Quad_Sa_Full}). Although we note that the algorithm of \citet{sagratella_computing_2016} ran on a less performing machine, the results in \cref{tab:Quad_S} highlight the remarkable effectiveness of \emph{ZERO Regrets}. The speedup in the performance seems to be considerably larger than the improvement associated with different hardware and software specifications (i.e., our algorithm is $100$ times faster in terms of computing times). 

\begin{table}[!ht]
\setlength\tabcolsep{0.7em}
\centering
\resizebox{0.50\textwidth}{!}{%
\begin{tabular}{>{\bf}clllll}
\toprule
\textbf{$\bm{n}$} & 
\textbf{\#EI} & 
\textbf{\#EQs} &
\textbf{Time} & 
\textbf{Time-\nth{1}} &
\textbf{Tl} \\
\midrule
2 & 81.33  & 3.67  & 2.21 & 0.58 & 0/12 \\
3 & 115.13 & 4.00  & 2.44 & 0.60 & 0/8  \\
4 & 119.00 & 10.25 & 3.90 & 0.97 & 0/4  \\
6 & 79.50  & 3.50  & 0.96 & 0.38 & 0/4   \\
\bottomrule \vspace{-\aboverulesep}
\end{tabular} %
}
\caption{Results for the \QIPG  from the instances of \citet{sagratella_computing_2016}. The complete table of results is in the Appendix (\cref{tab:Quad_Sa_Full}).}
\label{tab:Quad_S}
\end{table}

\section{Concluding Remarks}
\label{sec:concluding}
This paper presents a general framework to compute, enumerate and select equilibria for a class of \IPGs. These games are a fairly natural multi-agent extension of traditional problems in Operations Research, such as resource allocation, pricing, and combinatorial problems, and are powerful modeling tools for various applications. 
We provide a theoretical characterization of our framework through the concepts of \emph{equilibrium inequality} and \emph{equilibrium closure}. We explore the interplay between rationality and cutting planes by introducing a series of general and special-purpose classes of equilibrium inequalities and provide an interpretable criterion to frame the strategic interaction among players. The algorithm we introduce is general and it smoothly integrates with the existing optimization technology. Practically, we apply our framework to various problems from the relevant literature and significant application domains. We perform an extensive computational campaign and demonstrate the high efficiency and scalability of \emph{ZERO Regrets}. Our computational results also provide evidence of the existence of efficient \PNEs, further motivating the need for suitable algorithms to select or enumerate them. We also remark that our algorithm could practically work -- up to a numerical tolerance -- even when some of the players' variables are bounded and continuous, e.g., by dropping the integrality requirement on some variables.
We are prudently optimistic about the impact our framework may have in applied domains and the future methodological research directions it may open. We envision the potential for a series of theoretical contributions regarding the structure of new classes of general and problem-specific equilibrium inequalities and computational methods to further improve the algorithm's performance. Above all, we hope our framework will foster future academic research and clear the way for novel and impactful applications of \IPGs.

\iftoggle{ARXIV}{
    \subsubsection*{Acknowledgements}
We would like to thank Ulrich Pferschy for the valuable discussions concerning our work.
}{
    \ACKNOWLEDGMENT{We would like to thank Ulrich Pferschy for useful discussions on our work.}
}
\newpage
\bibliography{ZERO}
\label{sec:bib}

\newpage

\section*{Appendix}
\label{sec:Appendix}

\section*{\KPG Complexity Proof}
\label{sec:KPGProof}

\noindent We perform a reduction from the \emph{DeNegre Bilevel Knapsack Problem} ($BKP$) below, which is \SigmaTwoPC \citep{caprara_study_2014}. 
\begin{definition}[BKP]
Given two $m$-dimensional non-negative integer vectors $a$ and $b$ and two non-negative integers $A$ and $B$, the $BKP$ asks whether there exists a binary vector $x$ -- with $\sum_{j=1}^m a_jx_j \leq A$ -- satisfying $\sum_{j=1}^m b_jy_j(1-x_j) \leq B-1$ for any binary vector $y$ such that $\sum_{j=1}^m b_jy_j \leq B$.
\end{definition}

\noindent Without loss of generality, we assume $a_j \le A$ for any $j$. If this is not the case, we can always modify the original $BKP$ instance as follows:
\begin{enumerate*}
\item we replace $A$ with $2A+1$, any $a_j\le A$ with $2a_j$, and any $a_j>A$ with $(2A+1)$, and
\item we add a new element $m+1$ (i.e., new item), with $a_{m+1}=1$ and $b_{m+1}=B$. 
\end{enumerate*}
In any solution of this modified instance, we must have $x_{m+1}=1$, otherwise $\sum_{j=1}^{m+1} b_jy_j(1-x_j) \le B-1$ would never hold since $\sum_{j=1}^{m+1} b_jy_j(1-x_j) = B$ when $x_{m+1}=0$ and $y_{m+1}=1$. Setting  $x_{m+1}=1$ gives a residual capacity $2A$ for the packing constraint of $x$. Indeed, every subset of $x$ variables with original $a_j \le A$ that was satisfying $\sum_{j=1}^m a_jx_j \le A$ now satisfies $\sum_{j=1}^m 2a_jx_j \le 2A$. On the contrary, we cannot select any $x_j$ variable with original $a_j>A$. Thus, a solution (if any) to the original instance corresponds to a solution to the modified instance, and vice versa.
\iftoggle{ARXIV}{
    \begin{proof}
}{
    \proof{Proof of \cref{thm:IPGSigma2p}.}
}
First note that deciding if \KPG admits a \PNE is in $\Sigma^p_2$, as we ask whether there exists a strategy profile where every player cannot improve its payoff with any of its strategies, and we can compute the payoff of such strategies in polynomial time. Given a $BKP$ instance, we construct a \KPG instance with $2$ players as follows. We consider $m+1$ items and associate the elements of vectors $x$ and $y$ with the first $m$ elements of vectors $x^1$ and $x^2$, respectively. Then, player 1 solves the problem in (\ref{KPG:Proof:Player1}),
whereas player 2 solves the problem in (\ref{KPG:Proof:Player2}).

\begin{align}
\max_{x^1} \{ \sum_{j=1}^{m} b_{j}x_{j}^1x_{j}^2 + x_{m+1}^1x_{m+1}^2 : \sum_{j=1}^{m} a_{j}x_{j}^1 \le A, \; x^1\in\{0,1\}^{m+1}
\}.
\label{KPG:Proof:Player1}
\end{align}

\begin{align}
\max_{x^2} \{ (B - 1) x_{m+1}^2 + \sum\limits_{j=1}^{m} b_{j}x_{j}^2 - \sum_{j=1}^{m} b_{j}x_{j}^2x_{j}^1 : \nonumber \\ \sum_{j=1}^{m} b_{j}x_{j}^2 + Bx_{m+1}^2 \leq B, \; x^2 \in \{0,1\}^{m+1}
\}.
\label{KPG:Proof:Player2}
\end{align}

\noindent In order to prove the theorem, we show that the \KPG instance has a \PNE if and only if the corresponding $BKP$ instance admits a solution. 

\paragraph{$BKP$ admits a solution. } We assume the $BKP$ instance has a solution $\overline{x}$. We prove that $\hat{x }^1 = (\overline{x},1)$, $\hat{x}^2 = (\overline{0},1)$ (with $\overline{0}$ being an $m$-dimensional vector of zeros) is a \PNE. First, both the strategies $\hat{x}^1$ and $\hat{x}^2$ are feasible by construction. Given $\hat{x}^2$, player 1 attains the maximum payoff of $1$ by playing strategy $\hat{x}^1$. The strategy $\hat{x}^2$ yields a payoff of $B-1$ for player 2 when player 1 plays $\hat{x}^1$. Player 2 cannot profitably deviate by setting $x_{m+1}^2 = 0$. This follows from the fact that the $BKP$ instance has a solution $\overline{x}$ and, given that $\hat{x}^1_j = \overline{x}_j$ for $j=1,\dots, m$, the following inequality must hold.
\begin{align}
    \sum\limits_{j=1}^{m} b_{j}x_{j}^2 - \sum\limits_{j=1}^{m} b_{j}x_{j}^2\hat{x}_{j}^1 \leq B-1. \nonumber
\end{align}
Thus, the pair of strategies $(\hat{x}^1, \hat{x}^2)$ is also a \PNE for the \KPG instance.

\paragraph{$BKP$ has no solution. }
If the $BKP$ instance has no solution, player 2 never plays $x^2_{m+1} = 1$ in a best-response, as it can always obtain a payoff of $B$ with variables $x^2_1, \dots, x^2_m$ for any player 1's feasible strategy. Consider any player 2's best-response $\hat{x}^2$, with $\hat{x}^2_{m+1} = 0$, and assume the \KPG instance has a \PNE $(\hat{x}^1, \hat{x}^2)$. Then, in the player 1's best-response $\hat{x}^1$, there exists at least one $\hat{x}^1_j = 1$ when $\hat{x}_{j}^2 = 1$ and $b_j > 0$ (since $a_j \leq A$ for any $j$). However, in this case, player 2 would deviate from $\hat{x}^2$, since $\hat{x}^2$ gives a payoff $< B$ under $\hat{x}^1$. Thus, no \PNE exists in the \KPG instance.   
\iftoggle{ARXIV}{
    \end{proof}
}{
    \Halmos
    \endproof
} 

\section*{Extended Computational Results}
In the following sections, we report the full results for the our computational tests. The columns are similar to the ones reported in the previous tables, possibly with the following additions
\begin{enumerate*}
\item \emph{\#It} indicating the number of iterations of \emph{ZERO Regrets}, and
\item \emph{\PNE*} reporting the social welfare of the most efficient \PNE,
\item \emph{\PNE°} reporting the social welfare of the less efficient \PNE (if computed),
\item $OSW$ reporting the optimal social welfare in the game, and
\item \emph{Bound} reporting the last proven bound on $\mathcal{Q}$ before the latter becomes infeasible (or the algorithm hits the time limit), irrespective on whether the algorithm enumerated \PNEs or not.
\end{enumerate*}

\newpage

\subsection*{Full \KPG Results}
We report the two tables with the full \KPG results. In the first column of \cref{tab:KPG_Full2,tab:KPG_Full3} we add the field $I$ to specify the instance type. Specifically, the knapsack capacity of player $i$ is given by $\sum_{j=1}(w^i_j) I/10$.

\newpage

\begin{table}[H]
\setlength\tabcolsep{1em}
\centering
\resizebox{\textwidth}{!}{%
\begin{tabular}{>{\bf}crrrrrrrrrr}
\toprule
\textbf{($\bm n$, $\bm m$, $\bm d, \bm I$)}    & 
\textbf{\POS}  & 
\textbf{\#EI}  & 
\textbf{\#EI\_P} & 
\textbf{\#EI\_D} & 
\textbf{\#It} & 
\textbf{Time} & 
\textbf{Time-\nth{1}} & 
\textbf{\PNE*} & 
\textbf{$\bm{OSW}$} & 
\textbf{Bound}   \\
\midrule
(2, 25, A, 2)  & 1.106 & 12    & 0           & 8                  & 3          & 0.036    & 0.035       & 1884    & 2084  & 1884      \\
(2, 25, A, 5)  & 1.025 & 20    & 0           & 0                  & 3          & 0.095    & 0.093       & 3086    & 3163  & 3086      \\
(2, 25, A, 8)  & 1.000 & 12    & 0           & 1                  & 2          & 0.035    & 0.032       & 4883    & 4883  & 4883      \\
(2, 25, B, 2)  & 1.021 & 14    & 0           & 0                  & 3          & 0.067    & 0.065       & 1609    & 1643  & 1609      \\
(2, 25, B, 5)  & 1.025 & 28    & 0           & 4                  & 4          & 0.250    & 0.182       & 3456    & 3542  & 3459      \\
(2, 25, B, 8)  & -     & 10    & 0           & 7                  & 2          & 0.038    & 0.035       & 4624    & 4624  & 4624      \\
(2, 25, C, 2)  & -     & 16    & 7           & 5                  & 4          & 0.153    & -           & -       & 1480  & 1329      \\
(2, 25, C, 5)  & -     & 62    & 12          & 17                 & 11         & 0.967    & -           & -       & 2083  & 1863      \\
(2, 25, C, 8)  & 1.064 & 10    & 10          & 1                  & 4          & 0.037    & 0.036       & 2739    & 2914  & 2739      \\
(2, 50, A, 2)  & 1.024 & 24    & 0           & 3                  & 4          & 0.213    & 0.208       & 3824    & 3914  & 3824      \\
(2, 50, A, 5)  & 1.035 & 16    & 0           & 2                  & 3          & 0.214    & 0.212       & 6404    & 6626  & 6404      \\
(2, 50, A, 8)  & 1.016 & 20    & 0           & 3                  & 4          & 0.205    & 0.204       & 6703    & 6809  & 6703      \\
(2, 50, B, 2)  & 1.004 & 10    & 0           & 0                  & 2          & 0.043    & 0.040       & 3930    & 3946  & 3930      \\
(2, 50, B, 5)  & 1.004 & 42    & 0           & 34                 & 5          & 0.853    & 0.620       & 6931    & 6962  & 6936      \\
(2, 50, B, 8)  & 1.008 & 28    & 0           & 25                 & 6          & 0.620    & 0.501       & 9294    & 9372  & 9294      \\
(2, 50, C, 2)  & 1.018 & 8     & 25          & 0                  & 3          & 0.087    & 0.086       & 3173    & 3230  & 3173      \\
(2, 50, C, 5)  & -     & 112   & 25          & 22                 & 17         & 17.190   & -           & -       & 5654  & 4923      \\
(2, 50, C, 8)  & 1.134 & 98    & 31          & 12                 & 18         & 1.749    & 1.747       & 5358    & 6074  & 5358      \\
(2, 75, A, 2)  & 1.008 & 36    & 0           & 19                 & 4          & 0.407    & 0.401       & 5784    & 5831  & 5784      \\
(2, 75, A, 5)  & 1.004 & 40    & 0           & 49                 & 4          & 1.025    & 0.572       & 12701   & 12757 & 12702     \\
(2, 75, A, 8)  & 1.001 & 38    & 0           & 25                 & 3          & 0.359    & 0.342       & 16319   & 16337 & 16319     \\
(2, 75, B, 2)  & 1.033 & 122   & 0           & 41                 & 12         & 12.483   & 9.045       & 5690    & 5880  & 5694      \\
(2, 75, B, 5)  & 1.015 & 72    & 0           & 35                 & 8          & 5.865    & 1.420       & 10293   & 10449 & 10297     \\
(2, 75, B, 8)  & 1.010 & 108   & 0           & 26                 & 12         & 6.691    & 6.664       & 13769   & 13910 & 13769     \\
(2, 75, C, 2)  & 1.061 & 94    & 43          & 57                 & 9          & 3.072    & 3.064       & 4356    & 4623  & 4356      \\
(2, 75, C, 5)  & -     & 136   & 35          & 87                 & 18         & 134.899  & -           & -       & 7908  & 6934      \\
(2, 75, C, 8)  & 1.089 & 108   & 37          & 57                 & 18         & 5.289    & 4.865       & 8455    & 9207  & 8467      \\
(2, 100, A, 2) & 1.007 & 38    & 0           & 29                 & 5          & 1.409    & 1.153       & 8302    & 8357  & 8313      \\
(2, 100, A, 5) & 1.002 & 20    & 0           & 4                  & 2          & 0.355    & 0.188       & 18271   & 18301 & 18274     \\
(2, 100, A, 8) & 1.011 & 18    & 0           & 11                 & 3          & 0.521    & 0.398       & 18516   & 18723 & 18519     \\
(2, 100, B, 2) & 1.018 & 78    & 0           & 11                 & 8          & 4.294    & 4.281       & 8156    & 8303  & 8156      \\
(2, 100, B, 5) & 1.010 & 500   & 0           & 203                & 42         & 655.957  & 425.088     & 14246   & 14390 & 14248     \\
(2, 100, B, 8) & 1.002 & 38    & 0           & 25                 & 5          & 0.997    & 0.988       & 19054   & 19084 & 19054     \\
(2, 100, C, 2) & 1.048 & 116   & 49          & 33                 & 13         & 15.873   & 11.332      & 5808    & 6084  & 5817      \\
(2, 100, C, 5) & -     & 464   & 53          & 260                & 30         & tl & -           & -       & 9611  & 8958      \\
(2, 100, C, 8) & -     & 1512  & 64          & 66                 & 110        & tl & -           & -       & 9791  & 9007      \\
     \\
\bottomrule \\
\end{tabular}
}
\caption{Full results for the \KPG with $n=2$.}
\label{tab:KPG_Full2}
\end{table}

\newpage

\begin{table}[H]
\setlength\tabcolsep{1em}
\centering
\resizebox{\textwidth}{!}{%
\begin{tabular}{>{\bf}crrrrrrrrrr}
\toprule
\textbf{($\bm n$, $\bm m$, $\bm d, \bm I$)}    & 
\textbf{\POS}  & 
\textbf{\#EI}  & 
\textbf{\#EI\_P} & 
\textbf{\#EI\_D} & 
\textbf{\#It} & 
\textbf{Time} & 
\textbf{Time-\nth{1}} & 
\textbf{\PNE*} & 
\textbf{$\bm{OSW}$} & 
\textbf{Bound}   \\
\midrule
(3, 25, A, 2)  & 1.010 & 21    & 0           & 0                  & 3          & 0.166    & 0.164       & 3738    & 3777  & 3738      \\
(3, 25, A, 5)  & 1.004 & 30    & 0           & 0                  & 2          & 0.151    & 0.144       & 5480    & 5500  & 5480      \\
(3, 25, A, 8)  & 1.011 & 42    & 0           & 28                 & 3          & 0.323    & 0.316       & 9592    & 9693  & 9592      \\
(3, 25, B, 2)  & 1.034 & 27    & 0           & 3                  & 3          & 0.223    & 0.219       & 4535    & 4691  & 4535      \\
(3, 25, B, 5)  & 1.005 & 45    & 0           & 18                 & 3          & 0.394    & 0.388       & 7293    & 7329  & 7293      \\
(3, 25, B, 8)  & 1.008 & 60    & 0           & 23                 & 4          & 0.387    & 0.378       & 10346   & 10433 & 10346     \\
(3, 25, C, 2)  & 1.259 & 78    & 6           & 5                  & 8          & 6.765    & 5.643       & 2152    & 2710  & 2165      \\
(3, 25, C, 5)  & -     & 159   & 24          & 64                 & 13         & 82.115   & -           & -       & 4980  & 3935      \\
(3, 25, C, 8)  & -     & 36    & 59          & 32                 & 4          & 0.449    & -           & -       & 5735  & 4414      \\
(3, 50, A, 2)  & 1.033 & 99    & 0           & 17                 & 5          & 3.739    & 3.727       & 6769    & 6995  & 6769      \\
(3, 50, A, 5)  & 1.037 & 69    & 0           & 6                  & 5          & 2.413    & 2.043       & 11345   & 11764 & 11346     \\
(3, 50, A, 8)  & 1.007 & 117   & 0           & 50                 & 8          & 49.004   & 29.269      & 17283   & 17406 & 17283     \\
(3, 50, B, 2)  & 1.011 & 36    & 0           & 1                  & 4          & 1.976    & 1.971       & 7549    & 7634  & 7549      \\
(3, 50, B, 5)  & 1.015 & 468   & 0           & 99                 & 29         & tl & 483.842     & 13571   & 13781 & 13573     \\
(3, 50, B, 8)  & 1.011 & 114   & 0           & 33                 & 9          & 77.373   & 15.220      & 19680   & 19896 & 19697     \\
(3, 50, C, 2)  & -     & 231   & 37          & 108                & 15         & 934.599  & -           & -       & 5215  & 3764      \\
(3, 50, C, 5)  & -     & 159   & 64          & 397                & 10         & 211.139  & -           & -       & 9148  & 7564      \\
(3, 50, C, 8)  & -     & 54    & 88          & 169                & 5          & 0.977    & -           & -       & 11002 & 9194      \\
(3, 75, A, 2)  & 1.003 & 60    & 0           & 70                 & 4          & 9.057    & 1.342       & 14664   & 14711 & 14672     \\
(3, 75, A, 5)  & 1.041 & 45    & 0           & 1                  & 3          & 0.842    & 0.827       & 13869   & 14434 & 13869     \\
(3, 75, A, 8)  & 1.002 & 87    & 0           & 286                & 5          & 4.056    & 4.032       & 26468   & 26519 & 26468     \\
(3, 75, B, 2)  & -     & 444   & 0           & 130                & 18         & tl & -           & -       & 11508 & 11229     \\
(3, 75, B, 5)  & 1.002 & 81    & 0           & 97                 & 4          & 5.206    & 5.180       & 23139   & 23194 & 23139     \\
(3, 75, B, 8)  & 1.011 & 309   & 0           & 51                 & 18         & 1143.710 & 540.207     & 30118   & 30438 & 30118     \\
(3, 75, C, 2)  & -     & 357   & 36          & 177                & 15         & tl & -           & -       & 7242  & 6568      \\
(3, 75, C, 5)  & -     & 141   & 152         & 654                & 7          & 175.807  & -           & -       & 13553 & 11175     \\
(3, 75, C, 8)  & -     & 21    & 74          & 128                & 3          & 0.517    & -           & -       & 16736 & 14235     \\
(3, 100, A, 2) & -     & 333   & 0           & 15                 & 15         & tl & -           & -       & 15164 & 14825     \\
(3, 100, A, 5) & 1.003 & 408   & 0           & 28                 & 21         & tl & 1330.340    & 32673   & 32766 & 32677     \\
(3, 100, A, 8) & 1.006 & 42    & 0           & 391                & 3          & 1.959    & 1.915       & 37607   & 37826 & 37607     \\
(3, 100, B, 2) & -     & 516   & 0           & 297                & 21         & tl & -           & -       & 15946 & 15679     \\
(3, 100, B, 5) & -     & 291   & 0           & 81                 & 12         & tl & -           & -       & 29393 & 29119     \\
(3, 100, B, 8) & -     & 630   & 0           & 127                & 29         & tl & -           & -       & 40282 & 40082     \\
(3, 100, C, 2) & -     & 288   & 45          & 285                & 12         & tl & -           & -       & 11222 & 10045     \\
(3, 100, C, 5) & -     & 234   & 226         & 2059               & 10         & tl & -           & -       & 18272 & 15941     \\
(3, 100, C, 8) & -     & 30    & 242         & 715                & 3          & 0.932    & -           & -       & 20653 & 16855    
     \\
\bottomrule \\
\end{tabular}%
}
\caption{Full results for the \KPG with $n=3$.}
\label{tab:KPG_Full3}
\end{table}

\newpage

\subsection*{Full \NFG Results}
In \cref{tab:NFG_Full}, we report the full results for the \NFG. In the second and third columns, we report the weights of player $1$ and $2$ as $w^1$ and $w^2$, respectively. We remark that $w^3=1-w^1-w^2$.

\begin{table}[!ht]
\setlength\tabcolsep{1em}
\centering
\resizebox{\textwidth}{!}{%
\begin{tabular}{>{\bf}crrrrrrrrrr}
\textbf{($\bm{|V|}$, $\bm{|E|}$)}  & 
\textbf{$\bm{w^1}$}  & 
\textbf{$\bm{w^2}$} & 
\textbf{\POS}  & 
\textbf{\#EI}  & 
\textbf{\#It} & 
\textbf{Time} & 
\textbf{Time-\nth{1}} & 
\textbf{\PNE*} & 
\textbf{$\bm{OSW}$} & 
\textbf{Bound}   \\
\midrule
(50, 99)    & 0.33   & 0.33   & 1.061 & 5     & 3          & 0.037    & 0.036       & 980     & 924  & 980       \\
(50, 99)    & 0.6    & 0.2    & 1.245 & 8     & 3          & 0.040    & 0.039       & 1150    & 924  & 1150      \\
(50, 99)    & 0.45   & 0.45   & 1.061 & 5     & 3          & 0.034    & 0.034       & 980     & 924  & 980       \\
(100, 206)  & 0.33   & 0.33   & 1.000 & 3     & 2          & 0.047    & 0.041       & 1320    & 1320 & 1320      \\
(100, 206)  & 0.6    & 0.2    & 1.000 & 2     & 2          & 0.046    & 0.040       & 1320    & 1320 & 1320      \\
(100, 206)  & 0.45   & 0.45   & 1.000 & 2     & 2          & 0.047    & 0.040       & 1320    & 1320 & 1320      \\
(150, 308)  & 0.33   & 0.33   & 1.015 & 8     & 4          & 0.996    & 0.222       & 2049    & 2019 & 2042      \\
(150, 308)  & 0.6    & 0.2    & 1.015 & 5     & 4          & 0.354    & 0.353       & 2049    & 2019 & 2049      \\
(150, 308)  & 0.45   & 0.45   & 1.015 & 5     & 3          & 0.565    & 0.190       & 2049    & 2019 & 2041      \\
(200, 416)  & 0.33   & 0.33   & 1.000 & 1     & 2          & 0.109    & 0.096       & 2336    & 2336 & 2336      \\
(200, 416)  & 0.6    & 0.2    & 1.007 & 12    & 5          & 2.828    & 1.696       & 2352    & 2336 & 2346      \\
(200, 416)  & 0.45   & 0.45   & 1.187 & 22    & 10         & 6.908    & 1.529       & 2352    & 2336 & 2349      \\
(250, 517)  & 0.33   & 0.33   & 1.027 & 137   & 37         & 144.392  & 33.653      & 2730    & 2672 & 2729      \\
(250, 517)  & 0.6    & 0.2    & 1.027 & 47    & 17         & 43.991   & 13.430      & 2730    & 2672 & 2729      \\
(250, 517)  & 0.45   & 0.45   & 1.012 & 10    & 5          & 2.111    & 1.122       & 2703    & 2672 & 2693      \\
(300, 626)  & 0.33   & 0.33   & 1.060 & 36    & 10         & 14.877   & 2.068       & 3587    & 3567 & 3587      \\
(300, 626)  & 0.6    & 0.2    & 1.053 & 26    & 11         & 21.300   & 5.701       & 3587    & 3567 & 3585      \\
(300, 626)  & 0.45   & 0.45   & 1.000 & 1     & 2          & 0.161    & 0.140       & 3567    & 3567 & 3567      \\
(350, 730)  & 0.33   & 0.33   & 1.003 & 15    & 5          & 9.664    & 3.100       & 3678    & 3669 & 3677      \\
(350, 730)  & 0.6    & 0.2    & 1.014 & 41    & 11         & 31.889   & 18.997      & 3687    & 3669 & 3687      \\
(350, 730)  & 0.45   & 0.45   & 1.000 & 1     & 2          & 0.197    & 0.173       & 3669    & 3669 & 3669      \\
(400, 822)  & 0.33   & 0.33   & 1.207 & 100   & 29         & 163.047  & 0.228       & 4348    & 4319 & 4347      \\
(400, 822)  & 0.6    & 0.2    & 1.016 & 543   & 116        & tl & 584.854     & 4387    & 4319 & 4373      \\
(400, 822)  & 0.45   & 0.45   & 1.007 & 103   & 26         & 121.910  & 100.997     & 4348    & 4319 & 4346      \\
(450, 934)  & 0.33   & 0.33   & 1.159 & 0     & 2          & 0.304    & 0.250       & 4827    & 4827 & 4827      \\
(450, 934)  & 0.6    & 0.2    & 1.021 & 575   & 119        & tl & 7.284       & 4925    & 4827 & 4866      \\
(450, 934)  & 0.45   & 0.45   & 1.159 & 609   & 115        & tl & 0.281       & 4934    & 4827 & 4864      \\
(500, 1060) & 0.33   & 0.33   & 1.004 & 66    & 29         & 198.440  & 5.191       & 5535    & 5512 & 5534      \\
(500, 1060) & 0.6    & 0.2    & 1.004 & 20    & 8          & 20.951   & 11.231      & 5535    & 5512 & 5535      \\
(500, 1060) & 0.45   & 0.45   & 1.005 & 21    & 12         & 41.808   & 5.321       & 5535    & 5512 & 5534      \\
\bottomrule \\
\end{tabular} 
}
\caption{Full results for the \NFG.}
\label{tab:NFG_Full}
\end{table}

\newpage
\subsection*{Full \CFLD Results}

We report the results for a set of instances from \citet[Table 2]{cronert_equilibrium_2020} (i.e., $\beta=0.5$ and $d_{max}=20$). We report the full set of our results in \cref{tab:CFLD_Full}, where, in the second and third columns, we report the budget of player $1$ and $2$ as $B^1$ and $B^2$. When $n=3$, $B^3=10$.

{\footnotesize	
\setlength\tabcolsep{0.7em}
\begin{longtable}{>{\bf}c>{\bf}c>{\bf}ccrrrrrrrr}

\textbf{$\bm{n}$}  & 
\textbf{$\bm{B^1}$}  & 
\textbf{$\bm{B^2}$} & 
\textbf{\#EQs}  & 
\textbf{\POS}  & 
\textbf{\#EI}  & 
\textbf{\#It} & 
\textbf{Time} & 
\textbf{Time-\nth{1}} & 
\textbf{\PNE*} & 
\textbf{$\bm{OSW}$} & 
\textbf{Bound}   \\
\midrule
2 & 10 & 10  & 1     & 1.000 & 2     & 2          & 0.038    & 0.036       & 69      & 69  & 69        \\
2 & 10 & 20  & 1     & 1.000 & 2     & 2          & 0.259    & 0.256       & 109     & 109 & 109       \\
2 & 10 & 30  & 1     & 1.000 & 2     & 2          & 0.871    & 0.869       & 153     & 153 & 153       \\
2 & 10 & 40  & 1     & 1.000 & 2     & 2          & 1.026    & 1.025       & 186     & 186 & 186       \\
2 & 10 & 50  & 1     & 1.000 & 2     & 2          & 0.545    & 0.544       & 212     & 212 & 212       \\
2 & 10 & 60  & 1     & 1.000 & 2     & 2          & 0.627    & 0.626       & 232     & 232 & 232       \\
2 & 10 & 70  & 1     & 1.013 & 4     & 3          & 2.081    & 2.079       & 236     & 239 & 236       \\
2 & 10 & 80  & 1     & 1.047 & 8     & 4          & 4.908    & 4.905       & 236     & 247 & 236       \\
2 & 10 & 90  & 1     & 1.029 & 8     & 4          & 3.532    & 3.529       & 245     & 252 & 245       \\
2 & 10 & 100 & 1     & 1.028 & 8     & 4          & 3.706    & 3.701       & 247     & 254 & 247       \\
2 & 20 & 20  & 1     & 1.000 & 2     & 2          & 0.426    & 0.424       & 136     & 136 & 136       \\
2 & 20 & 30  & 1     & 1.000 & 2     & 2          & 1.153    & 1.151       & 180     & 180 & 180       \\
2 & 20 & 40  & 1     & 1.000 & 2     & 2          & 0.867    & 0.865       & 210     & 210 & 210       \\
2 & 20 & 50  & 1     & 1.000 & 4     & 3          & 1.760    & 1.758       & 232     & 232 & 232       \\
2 & 20 & 60  & 1     & 1.013 & 10    & 5          & 11.852   & 6.770       & 236     & 239 & 238       \\
2 & 20 & 70  & 1     & 1.038 & 6     & 3          & 7.494    & 7.194       & 234     & 243 & 234       \\
2 & 20 & 80  & 0     & -     & 6     & 4          & 9.530    & -           & -       & 252 & 243       \\
2 & 20 & 90  & 0     & -     & 8     & 4          & 14.304   & -           & -       & 254 & 247       \\
2 & 20 & 100 & 0     & -     & 6     & 3          & 19.163   & -           & -       & 254 & 252       \\
2 & 30 & 30  & 1     & 1.000 & 2     & 2          & 2.583    & 2.580       & 202     & 202 & 202       \\
2 & 30 & 40  & 1     & 1.000 & 2     & 2          & 1.852    & 1.849       & 232     & 232 & 232       \\
2 & 30 & 50  & 1     & 1.030 & 14    & 5          & 13.268   & 6.067       & 236     & 243 & 238       \\
2 & 30 & 60  & 1     & 1.065 & 14    & 7          & 37.077   & 37.067      & 232     & 247 & 232       \\
2 & 30 & 70  & 1     & 1.050 & 8     & 4          & 38.741   & 38.384      & 240     & 252 & 240       \\
2 & 30 & 80  & 1     & 1.058 & 16    & 5          & 515.179  & 270.395     & 240     & 254 & 241       \\
2 & 30 & 90  & 0     & -     & 10    & 6          & 1327.610 & -           & -       & 254 & 240       \\
2 & 30 & 100 & 0     & -     & 8     & 5          & 778.459  & -           & -       & 254 & 247       \\
2 & 40 & 40  & 1     & 1.138 & 24    & 9          & 491.695  & 154.949     & 210     & 239 & 216       \\
2 & 40 & 50  & 1     & 1.038 & 16    & 6          & 128.764  & 23.469      & 238     & 247 & 240       \\
2 & 40 & 60  & 2     & 1.050 & 18    & 9          & 344.539  & 98.475      & 240     & 252 & 240       \\
2 & 40 & 70  & 1     & 1.058 & 14    & 6          & 808.094  & 418.576     & 240     & 254 & 245       \\
2 & 40 & 80  & 1     & 1.058 & 10    & 5          & 1636.910 & 779.146     & 240     & 254 & 243       \\
3 & 10 & 10  & 1     & 1.072 & 6     & 3          & 0.360    & 0.358       & 69      & 74  & 69        \\
3 & 10 & 20  & 1     & 1.000 & 3     & 2          & 0.180    & 0.178       & 136     & 136 & 136       \\
3 & 10 & 30  & 1     & 1.000 & 3     & 2          & 0.522    & 0.518       & 180     & 180 & 180       \\
3 & 10 & 40  & 1     & 1.000 & 3     & 2          & 0.494    & 0.492       & 210     & 210 & 210       \\
3 & 10 & 50  & 1     & 1.000 & 3     & 2          & 0.631    & 0.628       & 232     & 232 & 232       \\
3 & 10 & 60  & 1     & 1.022 & 9     & 3          & 2.037    & 1.978       & 232     & 237 & 232       \\
3 & 10 & 70  & 1     & 1.030 & 9     & 4          & 4.772    & 4.769       & 236     & 243 & 236       \\
3 & 10 & 80  & 1     & 1.029 & 9     & 4          & 3.437    & 3.433       & 245     & 252 & 245       \\
3 & 10 & 90  & 1     & 1.037 & 9     & 4          & 6.679    & 6.676       & 245     & 254 & 245       \\
3 & 10 & 100 & 0     & -     & 9     & 4          & 16.522   & -           & -       & 254 & 249       \\
3 & 20 & 20  & 1     & 1.000 & 3     & 2          & 1.520    & 1.517       & 158     & 158 & 158       \\
3 & 20 & 30  & 2     & 1.037 & 18    & 4          & 6.161    & 5.795       & 187     & 194 & 187       \\
3 & 20 & 40  & 1     & 1.032 & 9     & 3          & 3.088    & 3.018       & 217     & 224 & 217       \\
3 & 20 & 50  & 1     & 1.000 & 3     & 2          & 1.931    & 1.929       & 239     & 239 & 239       \\
3 & 20 & 60  & 1     & 1.030 & 15    & 5          & 19.666   & 19.662      & 236     & 243 & 236       \\
3 & 20 & 70  & 1     & 1.068 & 6     & 3          & 5.114    & 5.111       & 236     & 252 & 236       \\
3 & 20 & 80  & 1     & 1.058 & 9     & 4          & 14.024   & 14.021      & 240     & 254 & 240       \\
3 & 20 & 90  & 0     & -     & 9     & 3          & 35.657   & -           & -       & 254 & 252       \\
3 & 20 & 100 & 0     & -     & 9     & 4          & 162.306  & -           & -       & 254 & 252       \\
3 & 30 & 30  & 1     & 1.000 & 18    & 4          & 11.841   & 11.838      & 216     & 216 & 216       \\
3 & 30 & 40  & 2     & 1.102 & 27    & 8          & 54.111   & 35.612      & 216     & 238 & 218       \\
3 & 30 & 50  & 1     & 1.038 & 24    & 8          & 51.141   & 51.135      & 238     & 247 & 238       \\
3 & 30 & 60  & 1     & 1.050 & 18    & 6          & 109.115  & 109.110     & 240     & 252 & 240       \\
3 & 30 & 70  & 1     & 1.058 & 24    & 8          & 226.185  & 226.180     & 240     & 254 & 240       \\
3 & 30 & 80  & 1     & 1.058 & 15    & 5          & 833.079  & 434.912     & 240     & 254 & 247       \\
3 & 30 & 90  & 0     & -     & 9     & 3          & 1153.500 & -           & -       & 254 & 254       \\
3 & 30 & 100 & 1     & 1.058 & 12    & 5          & tl & 3021.700    & 240     & 254 & 243       \\
3 & 40 & 40  & 1     & 1.038 & 36    & 8          & 222.003  & 221.639     & 238     & 247 & 238       \\
3 & 40 & 50  & 1     & 1.050 & 27    & 8          & 207.155  & 206.729     & 240     & 252 & 240       \\
3 & 40 & 60  & 1     & 1.076 & 24    & 9          & 2848.410 & 1696.300    & 236     & 254 & 244       \\
3 & 40 & 70  & 2     & 1.058 & 24    & 8          & 2441.450 & 833.853     & 240     & 254 & 240       \\
3 & 40 & 80  & 0     & -     & 15    & 6          & tl & -           & -       & 254 & 252          \\
\bottomrule \\
\caption{Full results for the \CFLD from the instances of \citet[Table 2]{cronert_equilibrium_2020}. \label{tab:CFLD_Full}}
\end{longtable}
}

\subsection*{Full \QIPGs Results}
We report the full results for the instances of \citet{schwarze_branch-and-prune_2022} in \cref{tab:Quad_SS_Full}, and the ones generated following the scheme of \citet{sagratella_computing_2016} in 
\cref{tab:Quad_Sa_Full}. In the latter table, we refer to \citet{sagratella_computing_2016} for an overview on the instances acronyms.

{\scriptsize
\setlength\tabcolsep{0.5em}
\begin{longtable}{>{\bf}ccrrrrrrrrrr}

\textbf{Instance}  & 
\textbf{\#EQs}  & 
\textbf{\POS}  & 
\textbf{\POA}  & 
\textbf{\#EI}  & 
\textbf{\#It} & 
\textbf{Time} & 
\textbf{Time-\nth{1}} & 
\textbf{\PNE*} & 
\textbf{\PNE°} & 
\textbf{$\bm{OSW}$} & 
\textbf{Bound}   \\
\midrule
C22\_3   & 2     & 1.0815  & 1.1238  & 14    & 4          & 0.3692  & 0.2835      & -22.7030  & -21.8475  & -24.5524  & -22.7030  \\
C22\_2   & 1     & 1.0000  & 1.0000  & 2     & 2          & 0.0949  & 0.0811      & -0.3146   & -0.3146   & -0.3146   & -0.3146   \\
C22\_1   & 2     & 1.4233  & 2.1559  & 24    & 6          & 0.6366  & 0.1852      & -13.5053  & -8.9158   & -19.2216  & 0.0196    \\
C22\_4   & 0     & -       & -       & 16    & 3          & 0.3097  & -           & -         & -         & -15.1462  & -6.8359   \\
C23\_1   & 2     & 1.0353  & 1.6333  & 28    & 5          & 0.7040  & 0.2652      & -10.7928  & -6.8413   & -11.1737  & 0.0000    \\
C23\_3   & 2     & 1.4506  & 3.0534  & 24    & 7          & 1.0306  & 0.5635      & -22.3566  & -10.6215  & -32.4315  & -10.6215  \\
C23\_2   & 0     & -       & -       & 26    & 6          & 1.0065  & -           & -         & -         & -22.3275  & -0.3407   \\
C23\_6   & 1     & 23.2815 & 23.2815 & 16    & 4          & 0.7469  & 0.6326      & -0.3396   & -0.3396   & -7.9063   & -0.3572   \\
C23\_7   & 1     & 1.0101  & 1.0101  & 6     & 3          & 0.4038  & 0.2817      & -4.5242   & -4.5242   & -4.5698   & -4.5242   \\
C23\_5   & 0     & -       & -       & 20    & 6          & 1.0198  & -           & -         & -         & -8.2644   & -0.2266   \\
C23\_4   & 0     & -       & -       & 26    & 6          & 0.9098  & -           & -         & -         & -44.4346  & -4.6428   \\
C23\_8   & 2     & 1.0153  & 1.5113  & 68    & 13         & 3.2209  & 0.3646      & -74.4543  & -50.0193  & -75.5936  & -3.5526   \\
C24\_4   & 0     & -       & -       & 32    & 5          & 4.6024  & -           & -         & -         & -46.2197  & -1.3690   \\
C24\_3   & 0     & -       & -       & 48    & 7          & 4.2759  & -           & -         & -         & -49.3061  & -2.2944   \\
C24\_2   & 0     & -       & -       & 40    & 5          & 2.8302  & -           & -         & -         & -50.0571  & -1.5728   \\
C24\_1   & 1     & 1.3206  & 1.3206  & 20    & 3          & 1.8845  & 1.1885      & -6.4656   & -6.4656   & -8.5384   & -6.4825   \\
C25\_4   & 1     & 1.4166  & 1.4166  & 34    & 7          & 34.4549 & 4.6820      & -50.3544  & -50.3544  & -71.3315  & -5.5704   \\
C25\_1   & 3     & 1.2068  & 11.3913 & 64    & 10         & 32.8601 & 3.6125      & -22.4829  & -2.3818   & -27.1321  & -2.3818   \\
C25\_3   & 1     & 2.0289  & 2.0289  & 60    & 10         & 26.9686 & 4.6411      & -45.6431  & -45.6431  & -92.6057  & -8.0414   \\
C25\_2   & 1     & 4.1130  & 4.1130  & 66    & 12         & 49.2170 & 17.9271     & -10.2162  & -10.2162  & -42.0192  & -2.1366   \\
C32\_1   & 2     & 3.1976  & 6.0787  & 30    & 4          & 0.7915  & 0.6078      & -21.6314  & -11.3788  & -69.1684  & -21.6314  \\
C32\_2   & 1     & 1.1101  & 1.1101  & 15    & 5          & 1.1502  & 0.4799      & -28.0541  & -28.0541  & -31.1421  & -9.9981   \\
C32\_3   & 3     & 1.1778  & -       & 63    & 9          & 2.2256  & 0.8736      & -45.5016  & 0.0000    & -53.5937  & -14.9792  \\
C32\_4   & 0     & -       & -       & 60    & 4          & 0.6875  & -           & -         & -         & -30.6557  & -8.8992   \\
C33\_2   & 1     & 8.3676  & 8.3676  & 117   & 12         & 23.1855 & 15.4598     & -9.2113   & -9.2113   & -77.0768  & -3.7804   \\
C33\_3   & 1     & 1.7099  & 1.7099  & 102   & 10         & 25.2901 & 1.6631      & -57.6349  & -57.6349  & -98.5491  & -1.3043   \\
C33\_1   & 1     & 1.2914  & 1.2914  & 129   & 17         & 55.9786 & 4.1433      & -138.1190 & -138.1190 & -178.3720 & -6.0334   \\
C33\_4   & 0     & -       & -       & 87    & 9          & 38.0066 & -           & -         & -         & -118.7970 & -6.0558   \\
NC22\_1  & 2     & 2.2947  & 2.2947  & 18    & 4          & 0.4126  & 0.1822      & -8.7456   & -8.7456   & -20.0687  & -8.7456   \\
NC22\_2  & 1     & 1.9081  & 1.9081  & 14    & 4          & 0.4483  & 0.3664      & -12.2614  & -12.2614  & -23.3957  & -12.2614  \\
NC22\_3  & 1     & 2.3939  & 2.3939  & 16    & 5          & 0.4243  & 0.3463      & -22.1224  & -22.1224  & -52.9584  & -22.1224  \\
NC22\_4  & 0     & -       & -       & 16    & 4          & 0.4330  & -           & -         & -         & -34.0944  & -22.9080  \\
NC23\_8  & 0     & -       & -       & 18    & 5          & 0.8855  & -           & -         & -         & -57.4117  & -31.8276  \\
NC23\_2  & 1     & 1.4346  & 1.4346  & 10    & 4          & 1.0951  & 0.2406      & -29.1437  & -29.1437  & -41.8083  & -15.2740  \\
NC23\_3  & 0     & -       & -       & 20    & 4          & 0.7570  & -           & -         & -         & -79.2272  & -28.3022  \\
NC23\_1  & 1     & 1.6194  & 1.6194  & 36    & 8          & 1.5463  & 0.6476      & -61.1489  & -61.1489  & -99.0215  & -2.1164   \\
NC23\_4  & 3     & 1.1508  & 1.9143  & 32    & 9          & 1.4491  & 0.3539      & -74.7629  & -44.9448  & -86.0367  & -35.2390  \\
NC23\_5  & 2     & 1.0962  & 1.7415  & 10    & 4          & 0.9394  & 0.3543      & -86.4907  & -54.4442  & -94.8133  & -54.4442  \\
NC23\_7  & 0     & -       & -       & 24    & 4          & 1.0328  & -           & -         & -         & -46.1839  & -14.7437  \\
NC23\_6  & 0     & -       & -       & 12    & 3          & 0.8541  & -           & -         & -         & -39.4816  & -29.1236  \\
NC24\_4  & 0     & -       & -       & 34    & 5          & 2.0809  & -           & -         & -         & -71.6710  & -62.7970  \\
NC24\_1  & 1     & 1.0190  & 1.0190  & 16    & 4          & 1.4601  & 0.8236      & -128.9180 & -128.9180 & -131.3660 & -98.6356  \\
NC24\_2  & 0     & -       & -       & 10    & 3          & 1.0317  & -           & -         & -         & -59.2505  & -50.3392  \\
NC24\_3  & 0     & -       & -       & 18    & 4          & 1.4827  & -           & -         & -         & -81.1047  & -62.3756  \\
NC25\_4  & 1     & 1.4370  & 1.4370  & 14    & 6          & 3.9086  & 2.9218      & -116.9060 & -116.9060 & -167.9990 & -116.9060 \\
NC25\_2  & 2     & 1.0487  & 1.1744  & 30    & 8          & 7.1719  & 0.9324      & -183.7380 & -164.0840 & -192.6940 & -72.1818  \\
NC25\_3  & 2     & 1.4358  & 1.9921  & 32    & 8          & 25.6082 & 3.2763      & -121.5220 & -87.5895  & -174.4870 & -87.5895  \\
NC25\_1  & 0     & -       & -       & 38    & 6          & 19.4117 & -           & -         & -         & -163.5600 & -62.2675  \\
NC32\_3  & 2     & 1.0000  & 1.4850  & 24    & 5          & 1.2589  & 0.4827      & -101.4570 & -68.3218  & -101.4570 & 2.6985    \\
NC32\_2  & 2     & 1.0652  & 1.4657  & 15    & 4          & 0.5213  & 0.2782      & -43.2125  & -31.4043  & -46.0281  & -31.4043  \\
NC32\_1  & 0     & -       & -       & 36    & 4          & 0.7526  & -           & -         & -         & -66.5208  & -21.4082  \\
NC32\_4  & 4     & 1.0145  & 1.8493  & 45    & 9          & 1.5802  & 0.3484      & -77.9484  & -42.7617  & -79.0771  & -20.4243  \\
NC33\_1  & 2     & 1.0042  & 1.0451  & 33    & 6          & 6.7043  & 1.3057      & -184.7260 & -177.4950 & -185.4930 & -99.5192  \\
NC33\_3  & 0     & -       & -       & 42    & 7          & 4.3455  & -           & -         & -         & -104.1130 & -21.0555  \\
NC33\_2  & 1     & 1.3586  & 1.3586  & 54    & 6          & 6.8289  & 3.3613      & -90.6533  & -90.6533  & -123.1600 & -71.4323  \\
NC33\_4  & 1     & 1.9431  & 1.9431  & 57    & 8          & 6.5654  & 2.6089      & -120.5760 & -120.5760 & -234.2880 & -41.7204     \\
\bottomrule \\
\caption{Full results for the \QIPG from the instances of \citet{schwarze_branch-and-prune_2022}. \label{tab:Quad_SS_Full}}
\end{longtable}
}

\begin{table}[!ht]
\setlength\tabcolsep{1em}
\centering
\resizebox{0.85\textwidth}{!}{%
\begin{tabular}{>{\bf}crrrrrrrrrrr}
\textbf{Instance}  & 
\textbf{\#EQs}  & 
\textbf{\#EI}  & 
\textbf{\#It} & 
\textbf{Time} & 
\textbf{Time-\nth{1}} & 
\textbf{\PNE*} & 
\textbf{\PNE°} & 
\textbf{$\bm{OSW}$} \\
\midrule
2-1-A-H  & 3     & 110   & 24         & 0.889  & 0.057       & -74.0   & 0.0      & -2128000.0 \\
2-1-B-H  & 2     & 102   & 16         & 0.581  & 0.158       & -8.5    & 0.0      & -8500000.0 \\
2-1-A-L  & 5     & 188   & 30         & 2.171  & 0.175       & -74.0   & 0.0      & -4376000.0 \\
2-1-B-L  & 2     & 136   & 23         & 0.743  & 0.100       & -27.5   & -0.5     & -6477570.0 \\
2-2-A-H  & 7     & 50    & 12         & 1.769  & 0.148       & -425.5  & 0.0      & -1438.0    \\
2-2-B-H  & 8     & 166   & 23         & 11.707 & 4.969       & -924.5  & -0.5     & -2712.0    \\
2-2-A-L  & 1     & 16    & 3          & 0.327  & 0.247       & 0.0     & 0.0      & -124.0     \\
2-2-B-L  & 7     & 112   & 19         & 6.500  & 0.205       & -7289.5 & 0.0      & -8560.0    \\
2-3-A-H  & 4     & 20    & 6          & 0.458  & 0.140       & -283.0  & 0.0      & -1118.0    \\
2-3-B-H  & 3     & 54    & 8          & 0.989  & 0.475       & -25.5   & 0.0      & -3138.0    \\
2-3-A-L  & 1     & 8     & 3          & 0.179  & 0.137       & 0.0     & 0.0      & -270.0     \\
2-3-B-L  & 1     & 14    & 3          & 0.250  & 0.188       & 0.0     & 0.0      & -750.5     \\
3-1-A-H  & 6     & 228   & 25         & 5.137  & 0.964       & -4776.0 & 0.0      & -76091.5   \\
3-1-B-H  & 8     & 246   & 31         & 4.465  & 0.647       & -957.0  & 0.0      & -234695.0  \\
3-1-A-L  & 3     & 159   & 18         & 3.220  & 1.273       & -618.5  & 0.0      & -93872.0   \\
3-1-B-L  & 1     & 105   & 12         & 1.204  & 0.392       & 0.0     & 0.0      & -71595.0   \\
3-2-A-H  & 1     & 33    & 5          & 0.760  & 0.395       & 0.0     & 0.0      & -1962.5    \\
3-2-B-H  & 1     & 15    & 3          & 0.078  & 0.069       & 0.0     & 0.0      & -1080.0    \\
3-2-A-L  & 8     & 84    & 11         & 3.390  & 0.629       & -1558.0 & 0.0      & -3032.5    \\
3-2-B-L  & 4     & 51    & 8          & 1.269  & 0.447       & -125.0  & 0.0      & -2044.0    \\
4-1-A-H  & 4     & 76    & 7          & 2.140  & 1.077       & -249.0  & 0.0      & -552.5     \\
4-1-B-H  & 13    & 152   & 16         & 4.654  & 1.689       & -3603.0 & 0.0      & -6115.5    \\
4-1-A-L  & 13    & 116   & 10         & 5.927  & 0.869       & -1462.0 & 0.0      & -1804.0    \\
4-1-B-L  & 11    & 132   & 12         & 2.863  & 0.238       & -1677.5 & 0.0      & -5817.5    \\
6-1-A-H  & 3     & 66    & 5          & 0.595  & 0.425       & -36.5   & 0.0      & -1437.5    \\
6-1-B-H  & 2     & 54    & 6          & 0.457  & 0.302       & -17.0   & 0.0      & -1715.0    \\
6-1-A-L  & 3     & 60    & 5          & 0.711  & 0.192       & -440.0  & 0.0      & -2795.0    \\
6-1-B-L  & 6     & 138   & 8          & 2.059  & 0.602       & -363.5  & 0.0      & -10510.0        \\
\bottomrule \\
\end{tabular} 
}
\caption{Full results for the \QIPG from the instances of \citet{sagratella_computing_2016}.}
\label{tab:Quad_Sa_Full}
\end{table}

\end{document}